%% file: fsewald.tex
\documentclass[a4paper,11pt]{article}

\usepackage{sectsty}
\sectionfont{\fontsize{12}{15}\selectfont}

\usepackage[font=small,labelfont=bf]{caption}
\usepackage{geometry}
\usepackage{amsmath,amssymb}
\usepackage{bm}
\usepackage{commath}
\usepackage{hyperref}
\usepackage[square,numbers]{natbib}
\usepackage{subcaption}
\usepackage{graphicx}

\hypersetup{
  pdfborder = 0 0 0,
  colorlinks=false,
  linkcolor=black,
  citecolor=black,
  filecolor=black,
  urlcolor=black,
} 

\input{defs}

\date{}

\title{Fast Ewald summation for free-space Stokes potentials}

\author{\small L. af Klinteberg\footnote{ludvigak@kth.se} , ~
  D. S. Shamshirgar\footnote{davoudss@kth.se} , ~
  A.-K. Tornberg\footnote{akto@kth.se}
\\~ \\
\small{\it
KTH Mathematics, Linn\'{e} Flow Centre/Swedish e-Science Research Centre, } \\
\small{\it100 44
 Stockholm, Sweden.}} 

\begin{document}

\maketitle

\noindent\makebox[\linewidth]{\rule{\textwidth}{1pt}} 
\section*{Abstract}
We present a spectrally accurate method for the rapid evaluation of free-space Stokes
potentials, i.e. sums involving a large number of free space Green's
functions. We consider sums involving stokeslets, stresslets and
rotlets that appear in boundary integral methods and potential methods for solving Stokes
equations.
The method combines the framework of the Spectral Ewald method for
periodic problems \cite{Lindbo2010}, with a very recent approach to solving the
free-space harmonic and biharmonic equations using fast Fourier
transforms (FFTs) on a uniform grid \cite{Vico2016}. 
Convolution with a truncated Gaussian function is used to place point
sources on a grid. 
With precomputation of a scalar grid quantity that does not depend on these sources, the amount of
oversampling of the grids with Gaussians can be kept at a factor of two,
the minimum for aperiodic convolutions by FFTs.
The resulting algorithm has a computational complexity of $O(N \log
N)$ for problems with $N$ sources and targets. 
Comparison is made with a fast multipole method (FMM) to show that the
performance of the new method is competitive. 

\noindent\makebox[\linewidth]{\rule{\textwidth}{1pt}} 

\setcounter{tocdepth}{2}

\section{Introduction}
 
In this paper, we consider the evaluation of free-space
potentials of Stokes flow, i.e. vector fields defined by sums involving a large
number of free space Green's functions such as the so-called  stokeslet, stresslet or rotlet. 
The stokeslet is the free space Green's function for velocity, and is given by 
\[
S(\rb)=\frac{1}{r} \id + \frac{1}{r^3} \rb \rb, \mbox{  or  }
S_{jl}(\rb)=\frac{\delta_{jl}}{r}  + \frac{r_j r_l}{r^3},  \quad j,l=1,2,3,
\]
with $r=|\rb|$ and where $\delta_{jl}$ is the Kronecker delta. 
The stresslet and rotlet will be introduced in the following. 
The discrete sums are on the form 
\begin{equation} 
\ub(\xb_{\mi})= 
\sum_{\substack{\ni=1 \\\ni \ne \mi}}^{N}
S(\xb_{\mi}-\xb_{\ni}) \fb(\xb_{\ni}), \quad \mi=1,\ldots,N.
\label{eqn:usumfree}
\end{equation}
and appear in boundary integral methods and potential methods for solving Stokes
equations.

These sums have the same structure as the classical Coulombic or
gravitational $N$-body problems that involve the harmonic kernel, and
the direct evaluation of such a sum for $\mi=1,\ldots,N$ requires
$O(N^2)$ work.  The Fast Multipole Method (FMM) can reduce that cost
to $O(N)$ work, where the constant multiplying $N$ will depend on the
required accuracy.  FMM was first introduced by Greengard and Rokhlin
for the harmonic kernel in 2D and later in 3D
\cite{Greengard1987,Cheng1999} and has since been extended to other
kernels, including the fundamental solutions of Stokes flow considered
here \cite{Rodin2000,Duraiswami2006,Tornberg2008,Wang2007,Ying2004}.  Related
is also the development of a so called pre-corrected FFT method based
on fast Fourier transforms. This method has been applied to the rapid
evaluation of stokeslet sums for panel-based discretizations of
surfaces \cite{White2006}.


For periodic problems, FFT-based fast methods built on the foundation
of so-called Ewald summation have been succesful. Also here,
development started for the harmonic potential, specifically for
evaluation of the electrostatic potential and force in connection to
molecular dynamic simulations, see e.g the survey by Deserno and Holm
\cite{Deserno1998}.  One early method was the {\em Particle Mesh Ewald
} (PME) method by Darden et al. \cite{Darden1993}, later refined to
the {\em Smooth Particle Mesh Ewald } (SPME) method by Essman et
al. \cite{Essmann1995}.  The SPME method was extended to the fast
evaluation of the stokeslet sum by Saintillan et
al. \cite{Saintillan2005}.  To recover the exponentially fast
convergence of the Ewald sums that is lost when such a traditional PME
approach is used, the present authors have developed a spectrally
accurate PME-type method, the Spectral Ewald (SE) method both for the
sum of stokeslets \cite{Lindbo2010}, and stresslets
\cite{AfKlinteberg2014a}.  It has also been implemented for the sum of
rotlets \cite{AfKlinteberg2016rot}, and the source code is available online \cite{se_github}.
The Spectral Ewald method was recently used to accelerate the
Stokesian Dynamics simulations in \cite{Wang2016}.


The present work deals with the efficient and fast summation of free
space Green's functions for Stokes flow (stokeslets, stresslets and
rotlets), as exemplified by the sum of stokeslets in
\eqref{eqn:usumfree}.  The problem has no periodicity, but the
approach will still be based on Ewald summation and fast Fourier
transforms (FFTs), using ideas from \cite{Vico2016} to extend the
Fourier treatment to the free-space case.  Before we
explain this further, we will introduce the idea behind Ewald
summation.

\subsection{Triply periodic Ewald summation}

Consider the Stokes equations in $\reals^3$, singularly forced at
arbitrary locations $\xb_{\ni}$, $\ni=1,\ldots,N$, with strengths $8
\pi \mu \fb(\xb_{\ni})\in \reals^3$ (with the $8 \pi \mu$ scaling for
convenience).  Introduce the three dimensional delta function
$\delta(\xb-\xb_0)$, and write
\begin{align}
-\nabla p + \mu \nabla^2 \ub + \gb(\xb) & =0, \quad 
\gb(\xb)= 8 \pi \mu \sum_{\ni=1}^N \fb(\xb_{\ni}) \, \delta(\xb-\xb_{\ni}),
\quad \quad
\quad \quad \quad \quad 
\label{eq:Stokes_sing} \\
\nabla \cdot \ub & = 0, 
\notag
\end{align}
where $\ub$ is the velocity, $p$ is the pressure and $\mu$ is the
viscosity. 
The free-space solution to this problem, evaluated at the source
locations, is given by \eqref{eqn:usumfree}. 

The classical Ewald summation formulas were derived for the triply
periodic problem for the electrostatic potential by Ewald in 1921
\cite{Ewald1921} and for the stokeslet by Hasimoto in 1959
\cite{Hasimoto1959}.  Here, assume that all the point forces are
located within a box $\domain=[-L_1 /2, L_1 /2] \times [-L_2 /2, L_2
/2] \times [-L_3 /2, L_3 /2]$, and that we impose periodic boundary
conditions.  The solution to this problem is a sum not only over all
the point forces, but also over all their periodic replicas,
\begin{equation} 
\ub^{3P}(\xb_{\mi})= \sum_{\pb \in P_3} \sum_{\ni=1}^{N*} 
S(\xb_{\mi}-\xb_{\ni}+\pb ) \fb(\xb_{\ni}), \quad \mi=1,\ldots,N.
\label{eqn:usumper}
\end{equation}
Here, the sum over $\pb$ formalizes the periodic replication of the
point forces with 
\begin{equation}  
P_3 =\{ (jL_1,lL_2,pL_3\}: (j,l,p) \in \Z^3 \}. \quad 
\label{eqn:perdef}
\end{equation}
The ${N*}$ indicates that the term ($\ni=\mi$, $\pb={\bf 0}$) is excluded
from the sum. 
The slow decay of the stokeslet however makes this infinite sum
divergent. 
To make sense of this summation, one usually assumes that the point
forces are balanced by a mean pressure gradient, such that the
velocity integrates to zero over the periodic box. 
Under these assumptions, Hasimoto \cite{Hasimoto1959} derived the following Ewald
summation formula
 \begin{align}
    \ub^{3P}(\xb_{\mi}) = &
    \sum_{\pb \in P_3}\ \sum_{\ni=1}^{N*}  S^R(\xb_{\mi} - \xb_{\ni} + \pb,\xi)  \fb(\xb_{\ni}) +
    \frac{1}{V} \sum_{\kb \neq 0}   \hat{S}^F(\kb,\xi)
    \sum_{\ni=1}^N \fb(\xb_{\ni}) \expk{(\xb_{\mi} - \xb_{\ni})}  
\notag \\
& +  \lim_{|\rb|\rightarrow 0} \left( S^{R}(\rb,\xi)-S(\rb) \right)\fb(\xb_{\mi}), 
\label{eq:usum_SR_SF_limit}
\end{align}
where the $\ni=\mi$, $\pb=0$ term is excluded from the real space
sum, $V=L_1 L_2 L_3$, and
\begin{equation}
  S^R(\rb, \xi)  = 2\left( \frac{\xi  e^{-\xi^2 r^2}  }{\sqrt{\pi} r^2} +
    \frac{ \erfc{(\xi r)}  }{2 r^3} \right) (r^2 \id + \rb\rb) -
  \frac{4\xi}{\sqrt{\pi}}  e^{-\xi^2 r^2}  \id, \quad , 
  \label{eq:SR_Hasimoto}
\end{equation}
\begin{equation}
  \hat{S}^F(\kb,\xi)=
  8\pi \left(1 + \frac{k^2}{4\xi^2} \right)\frac{1}{k^4}(\id k^2 - \kb
  \kb) e^{-k^2/4\xi^2} , 
  \label{eq:ShatH}
\end{equation}
with $r = | \rb|, \ k=|\kb|$,
\begin{align}
  \kb \in \mathbb K =  \{ 2\pi(j_1/L_1, j_2/L_2, j_3/L_3) : \v j \in \Z^3 \}
\end{align}
 and 
\begin{equation}
  \lim_{|\rb|\rightarrow 0} \left( S^{R}(\rb,\xi)-S(\rb)
  \right)=-\frac{4 \xi}{\sqrt{\pi}} \id.
\label{eq:Slimit}
\end{equation}
The last term  in \eqref{eq:usum_SR_SF_limit}  is commonly referred to as the self-interaction
term. When evaluating the potential at $\xb_{\mi}$, we should exclude
the contribution from the point force at that same location. 
For the real space part, we can directly
skip the term in the summation when $\pb=0$ and $\ni=\mi$.
We however need to subtract the contribution from this point that has
been included in the Fourier sum.  
We can use that $S^F=S-S^R$, and subtract the
limit as $|\rb| \rightarrow 0$   \eqref{eq:usum_SR_SF_limit}.  Both $S$ and $S^R$ are singular, but the limit
of the difference is finite \eqref{eq:Slimit}. 

Both sums now decay exponentially, one in real space and one in
Fourier space. The parameter $\xi>0$ is a decomposition parameter that controls the
decay of the terms in the two sums. The sum in real space can naturally be truncated to
exclude interactions that are now negligible. The sum in $k$-space
however, is still a sum of complexity $O(N^2)$, now with a very large
constant introduced by the sum over $\kb$.

Methods in the PME family make use of FFTs to evaluate the $k$-space
sum, accelerating the evaluation such that $\xi$ can be chosen
larger to push more work into the $k$-space sum, allowing for tighter
truncation of the real space sum, and in total an $O(N \log N)$
method. This procedure introduces approximations since a grid must be
used and, as with the FMM, the constant multiplying $N \log
N$ will depend on the accuracy requirements.

\subsection{The free-space problem and this contribution}

Considering the free space problem, we can introduce the same kind of
decomposition as in \eqref{eq:usum_SR_SF_limit}. 
The real space sum stays the same, with the minor
change that the sum over $\pb$ is removed, and the self interaction
term does not change. 
However, 
the discrete sum in Fourier space is replaced by the inverse
Fourier transform, 
\begin{align}
  \v u^F(\xb,\xi) = \frac{1}{(2\pi)^3} \int_{\mathbb R^3}
  \hat{\stokeslet}^F(\v k, \xi) 
\cdot \sum_{\ni=1}^N \fb(\xb_{\ni})  e^{i\v k\cdot (\xb-\v
    x_{\ni})} \dif \v k .
\label{eq:invtransf}
\end{align}
Here, note the $1/k^2$ singularity in $\hat{S}^F(\v k, \xi)$ as
defined in \eqref{eq:ShatH}. The integral is well defined, and
integration can be performed e.g. in spherical coordinates. A numerical
quadrature method in spherical coordinates would however require
non-uniform FFTs for non-rectangular grids in $k$-space.
Instead, we will use a very recent idea introduced by Vico et
al. \cite{Vico2016} to solve free space problems by FFTs on uniform
grids. 

The method by Vico et al. \cite{Vico2016} is based on the idea to use
a modified Green's function. With a right hand side of compact
support, and a given domain inside which the solution is to be found,
a truncated Green's function can be defined that coincides with the
original one for a large enough domain (and is zero elsewhere), such
that the analytical solution defined through a convolution of the
Green's function with the right hand side remains unchanged. 
The gain is that the Fourier transform of this truncated Green's
function will have a finite limit at $\kb=0$.  A length scale related
to the truncation will however be introduced, introducing oscillations
in Fourier space which will require some upsampling to resolve.

The authors of \cite{Vico2016} present this approach for radial Green's functions, e.g. 
the harmonic and biharmonic kernels. In the present work, we are
considering kernels that are not radial. We will however use this idea
in a substep of our method, defining the Fourier
transform of the biharmonic (for stokeslet and stresslet) or harmonic
(for rotlet) kernels, and define our non-radial kernels from these. 
The need of upsampling that the truncation brings can be taken care of
in a precomputation step, and hence for a scalar quantity only. What remains
is an aperiodic discrete convolution that requires an upsampling of a
factor of two.

The key ingredients in our method for the rapid summation of
kernels of Stokes flow (stokeslet, stresslet and rotlet) in free
space will hence be the following. 
We make use of the framework of Ewald summation, to split the
sums in two parts - one that decays rapidly in real space, and one
in Fourier space.  The Fourier space treatment is based on the
Spectral Ewald method for triply and doubly periodic problems that has
been developed previously \cite{AfKlinteberg2014a,Lindbo2010,
  Lindbo2011c,Lindbo2011e}.  This means that point forces will be
interpolated to a uniform grid using truncated Gaussian functions that
are scaled to allow for best possible accuracy given the size of the
support. The implementation of the gridding is made efficient by the
means of Fast Gaussian Gridding (FGG) \cite{Greengard2004,Lindbo2011c}. 

In the periodic problem, an FFT of each component of the grid function
is computed, a scaling is done in Fourier space, and after inverse
FFTs, truncated Gaussians are again used to evaluate the result at any
evaluation point.  The new development in this paper is to extend this
treatment to the free space case, when periodic sums are replaced by
discretized Fourier integrals. As mentioned above, a precomputation
will be made to compute a modified free-space harmonic or biharmonic
kernel that will be used to define the scaling in Fourier space.

The details are yet to be explained, but as we hope to convey in the
following, the method that we develop here for potentials of Stokes
flow, can easily be extended to other kernels. 
For any kernel that can be expressed as a differentiation of the harmonic and/or biharmonic
kernel, the Ewald summation formulas can easily be derived and only
minor changes in the implementation of the method will be needed. 

Any method based on Ewald summation and acceleration by FFTs will be
most efficient in the triply periodic case. As soon as there is one or
more directions that are not periodic, there will be a need of some
oversampling of FFTs, which will increase the computational cost.  For
the FMM, the opposite is true. The free space problem is the fastest
to compute, any periodicity will invoke an additional cost, which will
become substantial or even overwhelming if the base periodic box has a
large aspect ratio. Hence, implementing the FFT-based Spectral Ewald
method for a free-space problem and comparing it to an FMM method will
be the worst possible case for the SE method.  
Still, as we will show in the results section, using an open source
implementation of the FMM \cite{fmm_webpage}, our new 
method is competitive and often performs better than that implementation of the FMM for
uniform point distributions (one can however expect this adaptive FMM to
perform better for highly non-uniform distributions). 

There is an additional value in having a method that can be used for
different periodicities, thereby keeping the structure intact and
easing the integration with the rest of the simulation code,
concerning e.g. modifications of quadrature methods in a boundary
integral method to handle near interactions.  A three dimensional
adaptive FMM is also much more intricate to implement than the SE
method. Open source software for the Stokes FMM does exist for the
free space problem (as the one used here), but we are not
aware of any software for the periodic problem.


\subsection{Outline of paper}

The outline of the paper is as follows. In section
\ref{sect:freeGreen} we start by introducing the stokeslet, stresslet
and rotlet, and write them on the operator form that we will later
use. 
In section \ref{sect:Ewaldsum} we introduce the ideas behind
Ewald decomposition, and establish a framework for straight forward
derivation of decompositions of different kernels. 
The new approach to solve free-space problems by
FFTs introduced by Vico et al. \cite{Vico2016} is presented 
in the following section, together with a
detailed discussion on oversampling needs and precomputation.  The new
method for evaluating the Fourier space component is described in
section \ref{sect:kspace}, while the evaluation of the real space sum
is briefly commented on in section \ref{sec:eval-real-space}.
New truncation error estimates are derived in section
\ref{sec:truncation-errors}, and in section \ref{sect:summary} we
summarize the full method. Numerical results are presented in section
\ref{sect:numres}, where the performance of the method is discussed
and comparison to an open source implementation of the FMM
\cite{fmm_webpage} is made.

\section{Green's functions of free-space Stokes flow}
\label{sect:freeGreen}

We will consider three different Green's functions of free-space
Stokes flow, the stokeslet $\stokeslet$, the stresslet $\stresslet$
and the rotlet $\rotlet$. They are defined as
\begin{align}
  \stokeslet_{jl}(\v r) &= \frac{\delta_{jl}}{r} + \frac{r_jr_l}{r^3}, 
  \label{eq:def_stokeslet} \\
  \stresslet_{jlm}(\v r) &= -6 \frac{r_jr_lr_m}{r^5}
  \label{eq:def_stresslet}, \\
  \rotlet_{jl}(\v r) &= \epsilon_{jlm}\frac{r_m}{r^3},
  \label{eq:def_rotlet}
\end{align}
where $r = |\v r|$. They can equivalently be
formulated as operators acting on $r=|\rb|$ \cite{Pozrikidis1996,Fan1998},
and in the case of the rotlet also on $1/r$, 
\begin{align}
  \stokeslet_{jl} (\v r)&= \left(\delta_{jl}\nabla^2 - \nabla_j\nabla_l\right) r, 
  \label{eq:stokeslet_op}\\
  \stresslet_{jlm} (\v r)&= \left[ \left(
      \delta_{jl}\nabla_m+\delta_{lm}\nabla_j+\delta_{mj}\nabla_l
    \right) \nabla^2 - 2\nabla_j\nabla_l\nabla_m \right] r,
  \label{eq:stresslet_op}\\
  \rotlet_{jl}(\v r) &= \left(-\epsilon_{jlm}\nabla_m \nabla^2 \right)
  r =\left(-2 \epsilon_{jlm}\nabla_m \right) \frac{1}{r}. 
  \label{eq:rotlet_op}
\end{align}

For a single forcing term, $8\pi \mu \fb$ at a source location
$\xb_0$, we write the solution for velocity as 
\begin{equation}
\ub(\xb)= S(\xb-\xb_0) \fb, \mbox{ or }
u_j(\xb)= S_{jl}(\xb-\xb_0) f_l, \quad j=1,2,3,
\label{eq:vel_source_pt}
\end{equation}
where the repeated index is summed over according to Einstein summation convention.
Similarly to the velocity, the stress field and vorticity
associated with this velocity can be written, 
\[
\sigma_{jl}(\xb)= T_{jlm}(\xb-\xb_0) f_m, 
\quad
\omega_j(\xb)=\Omega_{jl}(\xb-\xb_0) f_l. 
\]
In integral equations, the stresslet often appears instead multiplying
sources with two indices, also producing a velocity, 
\[
\ub_{j}(\xb)= T_{jlm}(\xb-\xb_0) f_{lm}, 
\]
and this is the case that we will consider here. (The typical form is 
$T_{jlm} n_{l} q_m$, where ${\vec n}$ is a vector normal to a surface).  

We want to rapidly evaluate discrete-sum potentials of the type given
in (\ref{eqn:usumfree}), either at the source locations as indicated
in that sum, or at any other arbitrary points, and we want to do so
for the three different Green's functions.  To allow for a generic
notation in the following despite the differences, we introduce the
unconventional notation
\begin{align}
  \v u(\xb) = \sum_{\ni=1}^N \Green(\xb - \xb_{\ni}) \cdot \fb(\xb_{\ni}),
  \label{eq:usum}
\end{align}
where $\Green$ can denote either the stokeslet $\stokeslet$, the stresslet
$\stresslet$ and the rotlet $\rotlet$, and the dot-notation
$\ub=\Green(\rb) \cdot \fb$ will be understood to mean
\[
u_j(\xb)= \stokeslet_{jl}(\rb) f_l, \quad
u_j(\xb)= \stresslet_{jlm}(\rb) f_{lm}, \quad
u_j(\xb)= \rotlet_{jl}(\rb) f_l, \quad j=1,2,3, 
\]
in the three different cases.

\section{Ewald summation}
\label{sect:Ewaldsum}
\subsection{Decomposing the Green's function}
In Ewald summation we take a non-smooth and long-range Green's
function $\Green$, such as
(\ref{eq:def_stokeslet}--\ref{eq:def_rotlet}), and decompose it into two parts,
\begin{align}
   \Green(\v r) = \Green^R(\v r) + \Green^F(\v r).
\end{align}
This is done such that $\Green^R$, called the real space component,
decays exponentially in $r=|\rb|$. At the same time $\Green^F$, called Fourier
space component, decays exponentially in Fourier space. The original
example of this, derived by Ewald \cite{Ewald1921}, decomposes the
Laplace Green's function as
\begin{align}
  \frac{1}{r} = \frac{\erfc(\xi r)}{r} + \frac{\erf(\xi r)}{r},
  \label{eq:laplace_decomp}
\end{align}
where $\xi$ is a parameter that controls the decay rates in the real
and Fourier spaces. Here the real space component decays like
$e^{-\xi^2r^2}$, while the Fourier space component decays like
$k^{-2}e^{-k^2/4\xi^2}$. The rapid decay rates allow truncation of the
components; the real space component is reduced to local interactions
between near neighbors, while the Fourier space component is truncated
at some maximum wave number $\kmax$.

There are two different ways of deriving an Ewald decomposition, which
we shall refer to as screening and splitting. In screening, one
introduces a screening function $\gamma(\rb, \xi)$, $\int_{\mathbb
  R^3}\gamma(\rb, \xi) \dif \v r = 1$, that decays smoothly away from
zero. The Green's function is then decomposed using its convolution
with $\gamma$,
\begin{align}
  \Green(\v r) = \Green(\v r) - (\Green * \gamma)(\v r, \xi) + (\Green * \gamma)(\v r, \xi),
\end{align}
such that
\begin{align}
  \Green^R(\v r, \xi) &= \Green(\v r) - (\Green * \gamma)(\v r, \xi), \\
  \Green^F(\v r, \xi) &= (\Green * \gamma)(\v r, \xi),
\end{align}
and (by the convolution theorem)
\begin{align}
  \widehat{\Green}^F(\v k, \xi) &= \widehat{\Green}(\v k)\widehat\gamma(\v
  k, \xi), 
\label{eq:GF_Ggamma}
\end{align}
where $\widehat f$ denotes the Fourier transform of $f$,
\begin{align}
  \widehat f (\v k) = \ftrans[f](\v k) = \int_{\mathbb R^3} f(\xb)
  e^{-i \v k \cdot \xb} \dif \xb .
\end{align}

The original Ewald decomposition \eqref{eq:laplace_decomp} can be
derived in this fashion, using the screening function
\begin{align}
  \gamma_E(\v r,\xi)=\xi^3\pi^{-3/2}e^{-\xi^2r^2} \quad
  \rightleftharpoons \quad \widehat \gamma_E(\v k,\xi) =
  e^{-k^2/4\xi^2} , 
\end{align}
where $r=|\rb|$, $k=|\kb|$. 
For the stokeslet \eqref{eq:def_stokeslet} an Ewald decomposition was
derived by Hasimoto \cite{Hasimoto1959}, which was later shown
\cite{Hernandez-Ortiz2007} to be equivalent to using the screening
function
\begin{align}
  \gamma_H(\v r,\xi)=\xi^3\pi^{-3/2}e^{-\xi^2r^2} \left(\frac{5}{2}-\xi
    ^2 r^2\right) \quad \rightleftharpoons \quad \widehat
  \gamma_H(\v k,\xi) = e^{-k^2/4\xi^2}
  \left(1+\frac{1}{4}\frac{k^2}{\xi^2}\right) .
  \label{eq:hasimoto}
\end{align}
In splitting, one starts with the operator form of the Green's
function (\ref{eq:stokeslet_op}--\ref{eq:rotlet_op}), $\Green = \K r$,
and splits the Green's function using a splitting function
$\splitfun$,
\begin{align}
  \Green(\v r) = \K [r-\splitfun(r, \xi)] + \K\splitfun(r, \xi),
\end{align}
such that
\begin{align}
  \Green^R(\v r,\xi) & = \K[r-\splitfun(r,\xi)], \\
  \widehat{\Green}^F(\v k, \xi) &= \widehat\K(\v k)\widehat\splitfun(k,
  \xi),
\label{eq:GF_Khat_Phihat}
\end{align}
where $\widehat\K(\v k)$ denotes $\K$ applied to $e^{-i\v k \cdot \xb}$
(e.g. if $\K=\Delta$ then $\widehat\K=-|\v k|^2=-k^2$).  The splitting method
was invented by Beenakker \cite{Beenakker1986}, who used
\begin{align}
  \splitfun_B(r,\xi) = r \erf(\xi r) \quad \rightleftharpoons \quad
  \widehat\splitfun_B(k, \xi) = -\frac{8\pi}{k^4}\left( 1 +
    \frac{1}{4}\frac{k^2}{\xi^2} + \frac{1}{8}\frac{k^4}{\xi^4}
  \right) e^{-k^2/4\xi^2}.
\end{align}
We have now defined $ \widehat{\Green}^F(\v k, \xi)$ in two different
ways in (\ref{eq:GF_Ggamma}) and (\ref{eq:GF_Khat_Phihat}), and can
equate the two.  We have that $\Green=K|\rb|=Kr$, where
$B(\rb)=|\rb|=r$ is the fundamental solution of the biharmonic
equation, i.e.
\[
\nabla^4 \biharmonic(\rb)=- 8\pi \delta(\xb). 
\]
From this, we get 
\begin{equation}
\hat{\Green}(\kb) =\widehat{\biharmonic}(\kb)
\hat{K}(\kb)=-\frac{8\pi}{k^4}\hat{K}(\kb).
\label{eq:Ghat_wBhat}
\end{equation}
Hence, the screening and splitting methods can be shown
\cite{afKlinteberg2016phd} to be related to each other as
\begin{align}
  \widehat\splitfun = -\frac{8\pi}{k^4}\widehat\gamma
  \quad \rightleftharpoons \quad
  \gamma = -\frac{1}{8\pi}\nabla^4\splitfun .
\end{align}
Using this, one can derive the screening and splitting functions
related to the Ewald, Hasimoto and Beenakker decompositions, shown in
Table \ref{tab:decompositions}.
\begin{table}[htbp]
\resizebox{\textwidth}{!}{
  \centering
  \begin{tabular}{l|l|l|l}
      & $\gamma$ & $\widehat\gamma$ & $\splitfun$ \\
    \hline
    Ewald &
      $\alpha e^{-\xi^2r^2}$ 
    & $e^{-k^2/4\xi^2}$
    & $r\erf(\xi r) + \frac{r\erf(\xi r)}{2\xi^2r^2} + \frac{e^{-\xi^2r^2}}{\sqrt{\pi } \xi }$ \\
    Hasimoto &
      $\alpha e^{-\xi ^2r^2} \left(\frac{5}{2}-\xi ^2 r^2\right)$  
    & $e^{-k^2/4 \xi^2} \left(1+\frac{1}{4}\frac{k^2}{\xi^2}\right)$ 
    & $r\erf(\xi r) + \frac{e^{-\xi^2r^2}}{\sqrt{\pi } \xi }$ \\ 
    Beenakker &
      $ \alpha  e^{-\xi^2 r^2} \left(10-11 \xi ^2
      r^2+2 \xi ^4 r^4\right) $
    & $ e^{-k^2/4\xi^2}\left( 1 + \frac{1}{4}\frac{k^2}{\xi^2} +
      \frac{1}{8}\frac{k^4}{\xi^4} \right)$ 
    & $r\erf(\xi r)$     
  \end{tabular}}
  \caption{Summary of the screening and splitting functions related to the Ewald, Hasimoto and Beenakker decompositions (in that order). The constant $\alpha=\xi^3\pi^{-3/2}$, and $\xi$ is the decomposition parameter.} 
  \label{tab:decompositions}
\end{table}

The relations listed in the table are very useful in derivation of
Ewald summation formulas. Finding the real space part $\Green^R(\rb,\xi)$ 
 is easiest using the splitting approach, since this only involves differentiation. 
The $k$-space term is however simpler to derive with the screening
approach. Combining \eqref{eq:GF_Ggamma} and \eqref{eq:Ghat_wBhat}
it directly follows
\begin{align}
  \widehat{\Green}^F(\v k, \xi) &= \hat{K}(\kb) \widehat{\biharmonic}(\kb)
\widehat\gamma(\v  k, \xi). 
\label{eq:GF_KBGhat}
\end{align}
where $\widehat\K(\v k)$ denotes the prefactor that is produced as
$\K$ is applied to $e^{-i\v k \cdot \xb}$.

Considering the information in the table, we can see that the Ewald
decomposition yields the fastest decay in Fourier space. However, this
screening function can only be used if the Green's function can be
written as an operator acting on $1/r$, like the rotlet
\eqref{eq:rotlet_op}.
In this case, we can think about the splitting approach as if
splitting $1/r$ such that $\splitfun$ is $\erf(\xi r)/r$ as in
  \eqref{eq:laplace_decomp}. 

If we attempt to use the Ewald screening function for the stokeslet or
stresslet, this will not produce a useful decomposition since this
screening function does not ``screen'' the point forces. The field
produced by a point force convolved with the screening function does
not converge rapidly (with distance from the source location) to the
field produced by that point force. If we were to do the calculation,
this manifests itself in slowly decaying terms in the real space sum.

Both the Hasimoto and Beenacker screening functions work for the
stokeslet and stresslet. The Hasimoto decomposition will yield
somewhat faster decaying terms in both real and Fourier space, and
will henceforth be the one that we will use. 

\subsection{Ewald free-space formulas}
\label{sec:ewald-free-space}
In the triply periodic setting, the Ewald summation formula as derived
by Hasimoto was given in \eqref{eq:usum_SR_SF_limit}.  As given in
\eqref{eq:invtransf}, for the free-space problem the discrete sum in
Fourier space is replaced by the inverse Fourier transform.  With our
generic notation, we can evaluate the discrete-sum potential
\eqref{eq:usum} as
\begin{align}
  \v u(\xb) = \sum_{\ni=1}^N \Green^R(\xb - \xb_{\ni},\xi) \cdot \fb(\xb_{\ni}) + 
  \frac{1}{(2\pi)^3} \int_{\mathbb R^3} \widehat{
    \Green}^F(\v k,\xi) \cdot \sum_{n=1}^N \fb(\xb_{\ni}) e^{i\v k\cdot (\xb-\v
    x_{\ni})} \dif \v k .
  \label{eq:usumEwald}
\end{align}

We now apply the screening approach to derive the real space formulas,
\begin{align}
  \stokeslet_{jl}^R(\rb, \xi)  & = \left(\delta_{jl}\nabla^2 -
    \nabla_j\nabla_l\right) \left[ r-\splitfun^H(r,\xi) \right] 
\notag \\
 \stresslet^R_{jlm} (\v r) & = \left[ \left(
      \delta_{jl}\nabla_m+\delta_{lm}\nabla_j+\delta_{mj}\nabla_l
    \right) \nabla^2 - 2\nabla_j\nabla_l\nabla_m \right] \left[
    r-\splitfun^H(r,\xi) \right] 
\notag \\
  \rotlet^R_{jl}(\v r) & = - \epsilon_{jlm}\nabla_m \nabla^2 \left[ r-\splitfun^E(r,\xi) \right]=
-2 \epsilon_{jlm}\nabla_m 
  \frac{\erfc(\xi \rb)}{r} 
\notag 
\end{align}
where the splitting functions $\splitfun^H$ and $\splitfun^E$ are
found in the first and second lines of Table \ref{tab:decompositions}, 
and we get
\begin{align}
  \stokeslet_{jl}^R(\rb, \xi)  & = 2\left( \frac{\xi  e^{-\xi^2 r^2}  }{\sqrt{\pi}} +
      \frac{ \erfc{(\xi r)}  }{2 r} \right) (\delta_{jl} + \hat{r}_j \hat{r}_l) -
    \frac{4\xi}{\sqrt{\pi}}  e^{-\xi^2 r^2}  \delta_{jl} , \quad , 
\label{eq:SR_Op} \\
\stresslet^R_{jlm} (\v r) & = - \frac{2}{r} \left[ \frac{3 \, \erfc(\xi r)}{r} +
\frac{2 \xi}{\sqrt{\pi}}  \left(  3+2\xi^2r^2\right)
  e^{-\xi^2 r^2} \right]  \hat{r}_j \hat{r}_l  \hat{r}_m 
\notag \\
& \qquad \qquad \qquad \qquad +
\frac{4 \xi^3}{\sqrt{\pi}} e^{-\xi^2 r^2} (\delta_{jl}\hat{r}_m + \delta_{lm} \hat{r}_j
+\delta_{mj}\hat{r}_l), 
\label{eq:TR_Op} \\
  \rotlet^R_{jl}(\v r) & = 2 \varepsilon_{jlm} \hat{r}_m
\left( \frac{\erfc(\xi r)}{r^2} + \frac{2 \xi}{\sqrt{\pi}}
  \frac{1}{r} e^{-\xi^2 r^2} \right),
\label{eq:RR_Op} 
\end{align}
where $\hat{\rb}=\rb/|\rb|$. 
Only the stokeslet has a non-zero limit as given in 
\eqref{eq:Slimit}, that must be included to remove the self-interaction
when evaluating at a source point location. 

Let us now turn to the Fourier space terms.  We will use the screening
approach, with the Hasimoto screening for the stokeslet and stresslet,
and the Ewald screening for the rotlet.  From \eqref{eq:GF_KBGhat},
for each Green's function using the specific form of the differential
operator in \eqref{eq:stokeslet_op}-\eqref{eq:rotlet_op} to define
$\widehat\K(\v k)$, with $\widehat\gamma_E(\v k, \xi)$ and
$\widehat\gamma_H(\v k, \xi)$ as given in the first and second line of
Table \ref{tab:decompositions}, we obtain
\begin{align*}\small
 \widehat{\stokeslet}^F(\v k,\xi)  = \Gwoexp^{\stokeslet}(\v k,\xi)
 e^{-k^2/4 \xi^2}, \quad 
 \widehat{\stresslet}^F(\v k,\xi)  = \Gwoexp^{\stresslet}(\v k,\xi)
 e^{-k^2/4 \xi^2}, \quad 
 \widehat{\rotlet}^F(\v k,\xi)  = \Gwoexp^{\rotlet}(\v k,\xi)  e^{-k^2/4 \xi^2},
\end{align*}
where
\begin{align}
  \Gwoexp^{\stokeslet}_{jl}(\v k,\xi) &= 
  -\left(k^2\delta_{jl} - k_jk_l
  \right) 
  \left(1+k^2/(4\xi^2)\right) \widehat \biharmonic(|\kb|), 
\label{eq:AS}\\ 
\Gwoexp^{\stresslet}_{jlm}(\kb,\xi) &=
 i \left[(k_m \delta_{jl} +k_j \delta_{lm} +k_l  \delta_{mj}) k^2   -2
  k_j k_l k_m \right]  
  \left(1+k^2/(4\xi^2)\right) \widehat \biharmonic(|\v k|) , 
\label{eq:AT}\\
\Gwoexp^{\rotlet}_{jl} (\kb,\xi) &=
-  i \varepsilon_{jlm} k_m k^2 \widehat \biharmonic(|\v k|) =
2  i \varepsilon_{jlm} k_m \widehat \harmonic(|\v k|). 
\label{eq:AR}
\end{align}
Here, $\widehat{\harmonic}(k)$ and $\widehat{\biharmonic}(k)$ are
the Fourier transforms of the harmonic Green's function
$H(r)=1/r$ and the biharmonic Green's function $B(r)=r$,
respectively, 
\begin{equation}
\widehat{\harmonic}(k)
=\frac{4\pi}{k^2}, \quad \quad 
\widehat{\biharmonic}(k)
=-\frac{8\pi}{k^4}.
\label{eq:def_harm_biharm}
\end{equation}

For a smooth compactly supported function $\widehat{\Green}^F(\v
k,\xi)$, the Fourier integral in \eqref{eq:usumEwald} can be
approximated to spectral accuracy with a trapezoidal rule in each
coordinate direction, allowing for the use of FFTs for the
evaluation. Inserting the definitions in \eqref{eq:def_harm_biharm}
into (\ref{eq:AS}--\ref{eq:AR}), we can however note that the Fourier
component has a singularity at $k=0$ for all three Green's functions.

We will introduce modified Green's functions for the harmonic
and biharmonic equations, that will still yield the exact same result
as the original ones in the solution domain, and where the Fourier
transforms of these functions have no singularity for $k=0$. 
The necessary ideas will be introduced in the next session, following
the recent work by Vico et al.  \cite{Vico2016}.

\section{Free-space solution of the harmonic and biharmonic equations}
\label{sec:free-space-solution}

Consider the Poisson equation
\[
-\Delta \varphi(\xb)=4 \pi f(\xb)
\] with free space boundary conditions.  The solution is given by
\begin{equation} \varphi(\xb)=\int_{\reals^3} \harmonic(|\xb-\yb|)
f(\yb) \, d\yb=\frac{1}{(2\pi)^3} \int_{\reals^3}
\widehat{\harmonic}(|\kb|) \hat{f}(\kb) e^{i \kb \cdot \xb},
\label{eq:conv_harmonic}
\end{equation} where $\harmonic(r)=1/r$ is the harmonic Green's
function and $\widehat{\harmonic}(k)=4\pi/k^2$ its Fourier
transform. Note that they are both radial.

Assume now that $f$ is compactly supported within a domain
$\domaintilde$; a box with sides $\v \Ltilde$,
\begin{align} 
  \domaintilde = \{\xb \mid x_i \in [0, \Ltilde_i]\,\},
\end{align}
and that we seek the solution $\varphi(\xb)$ for $\xb \in
\domaintilde$.  The largest point to point distance in the domain is
$|\v\Ltilde|$.  Let $\trunc\ge|\v\Ltilde|$.  Without changing the
solution we can then replace $\harmonic$ with a truncated version,
\begin{align} 
  \harmonic^\trunc(r) = \harmonic(r) \rect\left(\frac{r}{2\trunc}\right),
\end{align}
where
\begin{align}
  \rect(x) =
  \begin{cases}
    1 &\text{for } |x| \le 1/2,\\
    0 &\text{for } |x| > 1/2.
  \end{cases} 
\end{align}
The Fourier transform of this truncated Green's function is 
\cite{Vico2016}
\begin{align}
  \widehat \harmonic^\trunc(k) = 8 \pi 
  \left( \frac{\sin(\trunc k/2)}{k}\right)^2.
  \label{eq:HhatR}
\end{align}
This function has a well-defined limit at $k=0$,
\begin{align}
  \widehat \harmonic^\trunc(0) = \lim_{k \rightarrow 0} \widehat
  \harmonic^\trunc(k) = 2 \pi \trunc^2.
\end{align}

Similarly, to solve the biharmonic equation on the same size domain,
we can define $\biharmonic^\trunc(r) = \biharmonic(r)
\rect\left(\frac{r}{2 \trunc}\right)$, which has the Fourier transform
\cite{Vico2016}
\begin{align} 
  \widehat \biharmonic^\trunc(k) = 4 \pi \frac{
    (2-\trunc^2k^2)\cos(\trunc k) + 2 \trunc k \sin(\trunc k) - 2
  }{k^4} ,
\label{eq:BhatR}
\end{align} with the limit value
\begin{align} \widehat\biharmonic^\trunc(0) = \lim_{k \rightarrow 0}
  \widehat \biharmonic^\trunc(k) = \pi \trunc^4.
\end{align}

\subsection{Solving the harmonic and biharmonic equations using FFTs}
\label{sec:FFTsolve} We will now describe how to solve the harmonic
and biharmonic equations using FFTs. In doing so, we will for
simplicity of notation assume that the domain is a cube,
i.e. $\Ltilde_1=\Ltilde_2=\Ltilde_3=\Ltilde$.

Then do the following
\begin{enumerate}
\item Introduce a grid of size $\Mtilde^3$ with grid size
$h=\Ltilde/\Mtilde$, and evaluate $f(\xb)$ on that grid.
\item Define an oversampling factor $\osamp$, and zero-pad to do a 3D
FFT of size $(\osamp \Mtilde)^3$, defining $\hat{f}(\kb)$ for
\[ \kb=\frac{2\pi}{\Ltilde} \frac{1}{\osamp} (k_1,k_2,k_3), \quad
\quad k_i \in \left\{ -\frac{\osamp \Mtilde}{2},\ldots, \frac{\osamp
\Mtilde}{2}-1 \right\}
\]
\item Set $\trunc=\sqrt{3}\Ltilde$ and evaluate $\widehat \harmonic^\trunc(k)$
\eqref{eq:HhatR}, with $k=|\kb|$ for the set of $\kb$-vectors defined
above.
\item Multiply $\hat{f}(\kb)$ and $\widehat H^\trunc(k)$ for each $\kb$. Do
a 3D IFFT and truncate the result to keep the $\Mtilde^3$ values
defining the approximation of the solution $\varphi(\xb)$ on the grid.
\end{enumerate} To solve the biharmonic equation instead, replace
$\widehat H^\trunc(k)$ by $\widehat B^\trunc(k)$ as given in
\eqref{eq:BhatR}.

With this, we have for all $\xb_j$ in the grid computed the approximation 
\[
\varphi(\xb_j)
\approx \frac{(\Delta k)^3}{(2\pi)^3} \sum_{k_1,k_2,k_3=-\osamp \frac{\Mtilde}{2}}^{\osamp \frac{\Mtilde}{2}-1}
\widehat{\harmonic}^\trunc(|\kb|) \hat{f}(\kb) e^{i \kb \cdot \xb_j} 
\]
where $\Delta k=\frac{2\pi}{\Ltilde} \frac{1}{\osamp}$.

Note that we at no occasion explicitly multiply with a
prefactor, assuming there is a built-in scaling of $1/(\osamp \Mtilde)^3$ in
the 3D inverse FFT.  There should be a multiplication with $h^3$ in step 2, 
and with $(\Delta k/2\pi)^3$ above, but that cancels such that only
the built-in scaling remains. 

Since the convolution is aperiodic, we need to oversample by at least
a factor of two. In Vico et al. \cite{Vico2016} they advise that we
need an additional factor of two to resolve the oscillatory behavior
of the Fourier transform of the truncated kernel, which would yield
$\osamp=4$. It does however turn out that the need of oversampling is
less than this, as we will discuss in the next section. 
If we oversample sufficiently, the error will decay spectrally with
$\Mtilde$ given that the right hand side $f$ is smooth.

\subsection{Zero-padding/oversampling}
\label{sec:oversampling}

Consider the first integral in \eqref{eq:conv_harmonic}. With $f$
compactly supported on a cube with size $\Ltilde^3$, $H$ must be
defined on a cube with size $(2\Ltilde)^3$ to be able to compute the
convolution. $H^\trunc$, with $\trunc=|\tilde{\v L}|=\sqrt{3}\Ltilde$
coincides with $H$ inside the sphere of radius $\trunc$, the smallest
sphere with the cube inscribed.  When we use the FFT, we ``periodize''
the computations. We hence need to zero pad the data so that this
periodization interval is large enough to make sure that $H^\trunc$ is
not polluted within the cube of size $(2\Ltilde)^3$.

\begin{figure}[htbp]
  \centering
  \includegraphics{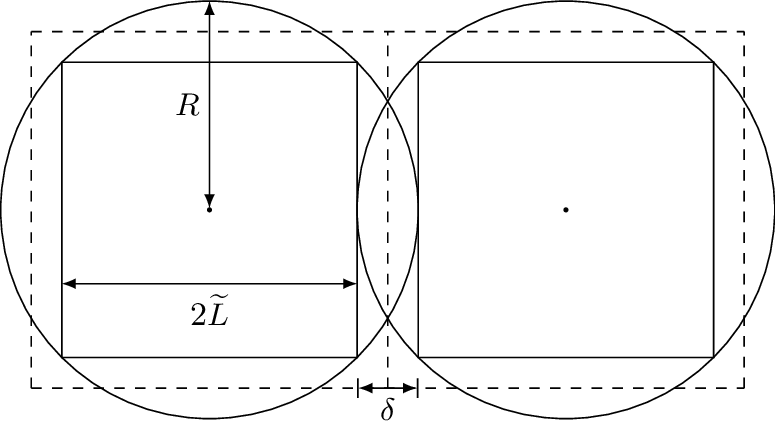}
  \caption{Illustration of the minimum zero-padding $\delta$ required
    to accurately represent the $\trunc$-truncated Green's function inside
    the domain of dimensions $2\widetilde{\v L}$, when using a
    periodic Fourier transform. The condition $\delta \ge \trunc-\Ltilde$ must be satisfied to avoid pollution from neighboring Green's
    functions inside the domain of interest.}
  \label{fig:oversampling}
\end{figure}

Assume that we zero pad the data up to a domain size of
$2\Ltilde+\delta$. If the sphere of radius $\trunc$ is to fit within
this domain, we would have $\delta=2(\trunc-\Ltilde)$. However, as is
illustrated in Figure \ref{fig:oversampling}, since it is enough that
$H^\trunc$ is not polluted within the cube, it is sufficient with
$\delta=\trunc-\Ltilde$.  In terms of an oversampling factor $\osamp$,
this corresponds to
\begin{align}
  \osamp \Ltilde  \ge 2 \Ltilde + \delta=
\Ltilde+\trunc, 
\end{align}
and the necessary condition becomes
\begin{align}
  \osamp \ge \frac{\Ltilde+ \trunc}{\Ltilde},
\end{align}
such that with $\trunc=\sqrt{3} \Ltilde$, we get 
$\osamp \ge 1 + \sqrt{3} \approx 2.8$.
Note, that an argument based instead on a large enough sampling ratio in the
Fourier domain to resolve the oscillatory truncated Green's function
would yield the smallest oversampling rate to be that with
$\delta=2(\trunc-\Ltilde)$, where the Green's function is without
pollution in the full sphere, and hence $\osamp \ge 2\sqrt{3} \approx 3.5$.

For non-cubic domains, we will have a larger oversampling requirement, 
\begin{align}
  \osamp \ge 1+ \frac{\trunc}{\min_i \Ltilde_i}, 
\end{align}
which  is $\osamp \ge 1 + |\widetilde{\v L}|/(\min_i \Ltilde_i)$ with the smallest possible $\trunc$. 
This additional cost can however be limited to a
precomputation step, through the scheme suggested in \cite{Vico2016}
as discussed in the next section. 

\subsection{Precomputation}
\label{sec:precomputation}

We will now further discuss step 4 in the algorithm introduced in
section \ref{sec:FFTsolve}. 
For ease of notation we do so in one
dimension. Each dimension will be treated the same way, so the
extension is simple to make. 
Let $\Mupgrid=\osamp \Mtilde$ be the number of grid points, such that 
$h=\Ltilde/\Mtilde=(\osamp \Ltilde)/\Mupgrid$, and let 
$k=(2\pi/(\osamp \Ltilde)) \kbar$, where $\kbar$ is an integer. 
By the means of an IFFT, we can compute
\[
\varphi_j=\frac{1}{\Mupgrid} \sum_{\kbar=-\Mupgrid/2}^{\Mupgrid/2-1}
\widehat{\Green}(k) \hat{f}_{\kbar} e^{i
  \frac{2\pi}{\Mupgrid} \kbar j}, \quad \quad j=0,\ldots,\Mupgrid-1,
\]
where $\widehat{\Green}(k)$ could be either
$\widehat{\harmonic}^\trunc(k)$ or $\widehat{\biharmonic}^\trunc(k)$,
and the Fourier coefficients $\hat{f}_{\kbar}$ have been computed by
an FFT,
\[
\hat{f}_{\kbar}=\sum_{l=0}^{\Mupgrid-1}
f(lh) e^{-i \frac{2\pi}{\Mupgrid} \kbar l} .
\]
Inserting how the $\hat{f}_{\kbar}$:s are actually computed, and
rearranging the order of the sums, we get
\[
\varphi_j=\sum_{l=0}^{\Mupgrid-1} 
\left[ \frac{1}{\Mupgrid} \sum_{\kbar=-\Mupgrid/2}^{\Mupgrid/2-1}
\widehat{\Green}(k)  e^{i  \frac{2\pi}{\Mupgrid} \kbar  (j-l)} \right]
f(lh) 
=\sum_{l=0}^{\Mupgrid-1} G_{j-l} f(lh), 
\]
where $G_{j-l}$, $j=0,\ldots,\Mupgrid-1$ will define the effective Green's function
on the grid, centered at grid point $l$. 
Note here, that $f$ has compact support and $f(lh)=0$ for
$l>\Mtilde-1$, and even though $\varphi_j$ is computed on the large grid, we
will truncate and keep only the $\Mtilde$ first values. 
Hence, for each $l$, only $\Mtilde$ values of $G_{j-l}$ is actually
needed to produce our result, and since $G_{(j+1)-(l+1)}=G_{j-l}$ a total
of $2\Mtilde$ grid values of the Green's function are used in the
calculation. 

\begin{figure}[htbp]
  \centering
  \begin{subfigure}[t]{0.49\textwidth}
    \centering
    \includegraphics[width=\textwidth,height=0.6\textwidth]{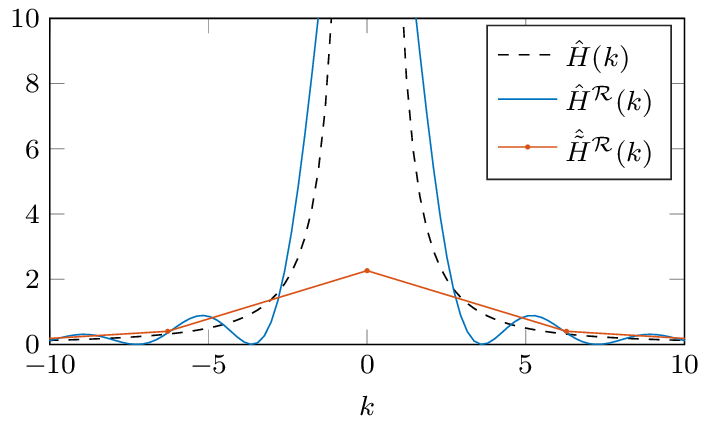}
    \caption{Fourier space representations of the original
      \eqref{eq:def_harm_biharm}, truncated \eqref{eq:HhatR} and
      mollified harmonic Green's functions. The latter is computed
      using an FFT.}
  \end{subfigure}
  \hfill
  \begin{subfigure}[t]{0.49\textwidth}
    \centering
    \includegraphics[width=\textwidth,height=0.6\textwidth]{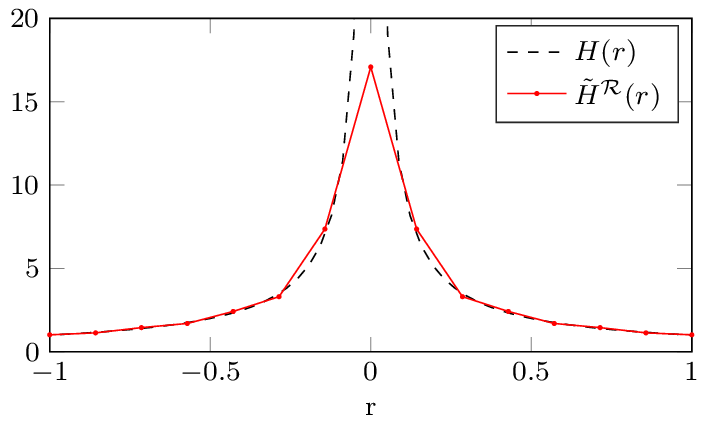}
    \caption{The harmonic Green's function and its mollified counterpart.}
  \end{subfigure}  
  \caption{Example of the mollified harmonic Green's function
    generated by inverse transform of $\hat{\harmonic}^\trunc$ using a
    finite-size IFFT.}
  \label{fig:mollified}
\end{figure}

Hence, it is possible to precompute an effective Green's function on
the grid , without knowing $f$, and truncate to save it centered on a
grid with $2\Mtilde$ values.  Let us denote by $\tilde{G}$ the
mollified Green's function that is the result of this procedure,
defined on a grid of size $2\Mtilde$, or in 3D $(2\Mtilde)^3$. An
example of the mollified harmonic Green's function is shown in Figure
\ref{fig:mollified}.

Once $f$ is given, we can now compute
our $\varphi_j$ values using an aperiodic convolution. 
This requires the following steps:
\begin{enumerate}
\item Do a 3D FFT of $\tilde{G}$.
\item Zero-pad and do a 3D FFT of size $(2\Mtilde)^3$ of $f$. 
\item Multiply the two, do an IFFT, and save the last $\Mtilde$ values
  in each dimension, i.e. $(\Mtilde)^3$ values as the result. 
\end{enumerate}

In practice, the FFT of $\tilde{G}$ is only done once if the equation
is to be solved with several right hand sides $f$.  When precomputing
$\tilde{G}$, an oversampling factor $\osamp \ge 1+\trunc/ \Ltilde$ is
needed as discussed in the previous section.  What remains is then an
aperiodic convolution with the right hand side as described above, and
an oversampling of $f$ with a factor of $2$ in this step is
sufficient.

\section{Evaluating the Fourier space component}
\label{sect:kspace}

Let us now go back to the Fourier space component in the Ewald
decomposition \eqref{eq:usumEwald}.  We will use the notation
\begin{align}
  \widehat{\Green}^{F,\trunc}(\v k, \xi) = \Gwoexp^{\Green,\trunc}(\v k, \xi)
  e^{-k^2/4\xi^2},
  \label{eq:Gwoexp}
\end{align}
where $\Green=\stokeslet$, $\stresslet$, and $\rotlet$, and were the
superindex $\trunc$ indicates that $\widehat{\harmonic}(k)$ and
$\widehat{\biharmonic}(k)$ are replaced by $\widehat{\harmonic}^\trunc(k)$
and $\widehat{\biharmonic}^\trunc(k)$ in the definitions \eqref{eq:AS},
\eqref{eq:AT} and \eqref{eq:AR}.  This means that the modified Green's
functions $\widehat{\stokeslet}^{F,\trunc}$, $\widehat{\stresslet}^{F,\trunc}$
and $\widehat{\rotlet}^{F,\trunc}$, will have no singularity at $k=0$.

The task is now to compute
\begin{align}
  \v u^F(\xb,\xi) = \frac{1}{(2\pi)^3} \int_{\mathbb R^3} e^{i\v
    k\cdot\xb} e^{-k^2/4\xi^2}
  \Gwoexp^{\Green,\trunc}(\v k, \xi) \cdot \sum_{n=1}^N \fb(\xb_{\ni}) e^{-i\v
    k \cdot \xb_{\ni}} \dif \v k 
  \label{eq:uf_sum}
\end{align}
for a given set of target points. The integrand of the inverse
transform is now smooth and can after truncation be easily evaluated
using the trapezoidal rule, but evaluation is still costly ---
$\ordo(N^2)$ if evaluating at $N$ target points. We will now outline
the spectral Ewald method, which uses the fast Fourier transform to
reduce the cost of this evaluation, yielding a method with a total
cost (including the real space sum) of $\ordo(N \log N)$. Before we
discuss the actual discretization and implementation details, we start
by describing the mathematical foundation of the method.

\subsection{Foundations}

First we introduce a scalar parameter $\eta>0$ and define
\begin{align}
  \widehat{\gridfcn}(\v k, \xi) = \sum_{n=1}^N \fb(\xb_{\ni}) e^{-i\v k\cdot
    \xb_{\ni}} e^{-\eta k^2/8\xi^2},
\label{eq:gridfcnhat_def}
\end{align}
which is the Fourier transform of the smooth function
\begin{align}
  \gridfcn(\xb, \xi) = \sum_{n=1}^N \fb(\xb_{\ni}) \left(\frac{2\xi^2}{\pi\eta}\right)^{3/2}
  e^{-2\xi^2|\xb-\xb_{\ni}|^2/\eta}.
  \label{eq:Hx}
\end{align}
Hence, instead of using  \eqref{eq:gridfcnhat_def} to directly evaluate
$\widehat{\gridfcn}(\v k, \xi)$, we can evaluate  $\gridfcn(\xb, \xi)$
and use an actual computation of the Fourier transform, 
\begin{equation}
  \widehat{\gridfcn}(\v k, \xi) =\int_{\reals^3} \gridfcn(\xb, \xi)
  e^{-i \kb \cdot \xb} \, \dif \xb.
\end{equation}
We furthermore define 
\begin{align}
  \widehat{\gridfcnsc}(\v k, \xi) = e^{-(1-\eta)k^2/4\xi^2}
  \Gwoexp^{\Green,\trunc}(\v k, \xi) \cdot \widehat{\gridfcn}(\v k,
  \xi),
  \label{eq:gridfcnsc}
\end{align}
such that   \eqref{eq:uf_sum} can be written
\begin{align}
  \v u^F(\xb, \xi) = \frac{1}{(2\pi)^3} \int_{\mathbb R^3}
  \widehat{\gridfcnsc}(\v k, \xi) e^{-\eta k^2/8\xi^2}  e^{ i\v
    k\cdot\xb } \dif \v k .
  \label{eq:Fx}
\end{align}
Using the convolution theorem we can write this as
\begin{align}
  \v u^F(\xb, \xi) = \int_{\mathbb R^3} \gridfcnsc(\v y,\xi)
  \left(\frac{2\xi^2}{\pi\eta}\right)^{3/2} e^{-2\xi^2|\xb-\v
    y|^2/\eta} \dif \v y, 
  \label{eq:u_integr}
\end{align}
where
\begin{equation}
  \gridfcnsc(\xb,\xi)
  = \frac{1}{(2\pi)^3} \int_{\reals^3} \widehat{\gridfcnsc}(\v k, \xi)
  e^{i \kb \cdot \xb} \, \dif \kb. 
\end{equation}

\subsection{Discretization}
\label{sec:discretization}

Assume that we are to evaluate \eqref{eq:uf_sum} for $\xb=\xb_{\mi}$,
$\mi=1\ldots N$, and for simplicity of notation that all points are contained in a
cube with equal sides $L$,
\[
  \xb_{\ni} \in \domain = [0, L]^3, \quad \ni = 1 \dots N. 
\]
The choice of $\eta$ will be discussed shortly, in section
\ref{sec:approximation-errors}. At this point, assume that the
Gaussians $e^{-2\xi^2|\cdot|^2/\eta}$ in \eqref{eq:Hx} and
\eqref{eq:u_integr} decay rapidly, and will be truncated outside a
diameter $\supp$. Then $\gridfcn$ becomes compactly supported, such
that we can compute $\widehat{\gridfcn}$ using an FFT, and the
integral in \eqref{eq:u_integr} becomes a local operation around each
target point $\xb_{\mi}$. To accomodate the support of the truncated
Gaussians, we must extend the domain by some length $\deltaL$. We will
discuss the choice of this length in the discussion on $\eta$. For
now, we consider the extended domain with sides $\Ltilde = L +
\deltaL$,
\begin{align}
  \domaintilde = [-\deltaL/2, L + \deltaL/2]^3.
\end{align}
This domain is discretized using a uniform grid with $\Mtilde^3$
points and grid spacing $h = \Ltilde / \Mtilde$.

To initialize our calculations, we precompute
$\widehat{\harmonic}^\trunc(k)$ in case of the rotlet, and
$\widehat{\biharmonic}^\trunc(k)$ in case of the stokeslet or
stresslet, as described in section \ref{sec:precomputation}. They need
to be precomputed on a domain of size $2\Ltilde$, with
$\trunc = \sqrt{3}\Ltilde$.


The first step of our computations is to evaluate $\gridfcn$ on the
grid as in \eqref{eq:Hx}. After that we zero-pad the FFT by a a
factor of $2$, to have an oversampled representation of
$\widehat{\gridfcn}$, before we scale it to define
$\widehat{\gridfcnsc}$ as in \eqref{eq:gridfcnsc}.  We will then
multiply by the precomputed fundamental solution
($\widehat{\harmonic}^\trunc(k)$ or $\widehat{\biharmonic}^\trunc(k)$)
and the additional scaling factors as given in \eqref{eq:AS},
\eqref{eq:AT} and \eqref{eq:AR}, and apply an inverse FFT to perform a
discrete convolution.

The computation of $\v u^F(\xb_{\mi},\xi)$, $\mi=1,\ldots,N$ can hence
be broken down into the following steps:
\begin{enumerate}
\item \emph{Spreading}: Compute $\gridfcn$ on the grid using
  \eqref{eq:Hx} and truncated Gaussians.
\item \emph{FFT}: Compute $\widehat{\gridfcn}$ using the
  three-dimensional FFT, zero-padded to the double size.
\item \emph{Scaling}: Compute $\widehat{\gridfcnsc}$ using
  \eqref{eq:gridfcnsc} and precomputed $\widehat{\harmonic}^\trunc(k)$
  or $\widehat{\biharmonic}^\trunc(k)$.
\item \emph{IFFT}: Apply the inverse three-dimensional FFT to
  $\widehat{\gridfcnsc}$. Truncate the result to have $\gridfcnsc$
  defined on the original grid.
\item \emph{Quadrature}: For each $\xb_{\mi}$, $\mi=1,\ldots,N$
  evaluate $\v u^F$ using \eqref{eq:u_integr} and the trapezoidal rule,
  with the Gaussian truncated outside the sphere of diameter $\supp$ centered at
  $\xb_{\mi}$.
\end{enumerate} This is the spectral Ewald method. A major cost of the
method is the large number of exponential function evaluations in
steps 1 and 5. This can be accelerated through the method of fast
Gaussian gridding (FGG) \cite{Greengard2004,Lindbo2011c}. It is then
natural to truncate the Gaussians outside a cube of $P^3$ grid points,
in which case
\begin{align}
  \supp = hP .
\end{align}
The computational cost of the FGG in steps 1 and 5 is then $\mathcal
O(NP^3)$, while the cost of the FFTs in steps 2 and 4 is
$\ordo(\Mtilde^3 \log \Mtilde )$. The cost of the scaling in step 3 is
$\ordo(\Mtilde^3)$, and negligible in this context.

\subsection{Errors in the spectral Ewald method}
\label{sec:approximation-errors}

The use of the spectral Ewald method for computing the Fourier space
component introduces approximation errors in the solution, which are
separate from the Fourier integral truncation error (further discussed
in section \ref{sec:fourier-space-error}). The approximation errors
stem from the use of a discrete quadrature rule in the quadrature
step, and from the truncation and discretization of the Gaussians
$e^{-2\xi^2|\cdot|^2 / \eta}$ in the spreading and quadrature
steps. The Ewald parameter $\xi$ should be regarded as free, since it
is used for work and error balancing between the real and Fourier
space sums (more on this in section \ref{sec:parameter-choice}). This
leaves two variables for controlling the approximation errors: the
scalar parameter $\eta$ and the Gaussian truncation width
$\supp=hP$. Following \cite{Lindbo2011c}, we write $\eta$ as
\begin{align}
  \eta = \left( \frac{\xi d}{m} \right)^2,
\end{align}
where $m$ is a shape parameter controlling how fast the Gaussian
decays within the support $\supp$. It can be shown \cite{Lindbo2011c}
that the approximation errors decay exponentially in $P$ with the
choice $m(P) = C \sqrt{\pi P}$, and that the constant $C$ should be
taken slightly below unity for optimal results (we use the value
$C=0.976$ suggested in \cite{Lindbo2011c}). With these choices, the
approximation errors of the method are controlled through a single
parameter $P$, and they furthermore decay exponentially in that
parameter.

\begin{figure}[htbp]
  \centering 
  \begin{subfigure}[b]{0.49\textwidth}
    \centering
    \includegraphics[width=\textwidth,height=0.6\textwidth]{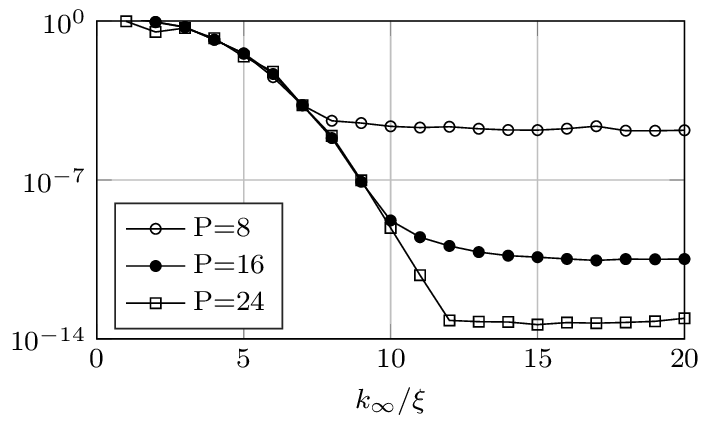}
    \caption{$\deltaL$ set through \eqref{eq:deltaL_cond}.}
  \end{subfigure}
  \hfill
  \begin{subfigure}[b]{0.49\textwidth}
    \centering
    \includegraphics[width=\textwidth,height=0.6\textwidth]{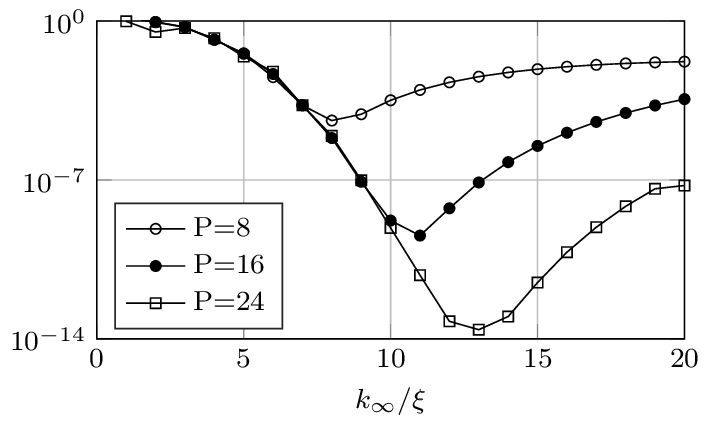}
    \caption{$\deltaL=d$.}
  \end{subfigure}
  \caption{Error in the stokeslet Fourier space component for various
    values of the discrete Gaussian support $P$. The system is a unit
    cube with 1000 sources, $\Mtilde \in [2,40]$, $\xi=2\pi$ and
    $\kmax = \pi \Mtilde / L=\pi \Mtilde$. To the left the domain is extended to
    include the support of the remainder Gaussians, while to the right
    the domain is only extended to cover the support of the gridding
    and quadrature Gaussians. Evidently, the extra support is only
    needed if $\Mtilde$ is picked larger than necessary for a given
    error tolerance.}
  \label{fig:P_errors}
\end{figure}

It is evident from the algorithm that $\deltaL$ must be chosen such
that the support of the truncated Gaussians in \eqref{eq:Hx} and
\eqref{eq:u_integr} is included, i.e. $\deltaL \ge \supp$. However, it
turns out that this is not always enough. In the spectral Ewald method
we have taken the Gaussian $e^{-\xi^2 r^2}$ of the screening function
(Table \ref{tab:decompositions}) and separated it into a series of
convolutions of Gaussians, through the factorization
\begin{align}
  e^{-k^2/4\xi^2}=e^{-\eta k^2/8\xi^2} \cdot e^{-(1-\eta)k^2/4\xi^2}
  \cdot e^{-\eta k^2/8\xi^2}  .
\end{align}
The first and last factors correspond to the Gaussian
$e^{-2\xi^2|\cdot|^2 / \eta}$ in the gridding and quadrature steps,
and is already properly resolved and truncated by our choices of
$\eta$ and $d$. For $\eta \ge 1$ the entire ''original'' Gaussian is
contained in these two factors, and the middle factor can be viewed as
a deconvolution of the type used in the nonuniform FFT
\cite{Lee2005}. However, for $\eta < 1$ the middle factor corresponds
to the ''remainder'' Gaussian $e^{-\xi^2|\cdot|^2 / (1 - \eta)}$, and
\eqref{eq:gridfcnsc} corresponds to a convolution with that Gaussian,
carried out in Fourier space. For the convolution to be properly
represented, we must make sure that the domain $\Ltilde$ includes the
support of $e^{-\xi^2|\cdot|^2 / (1 - \eta)}$ to the desired
truncation level. The original Gaussians are truncated at the level
$e^{-2\xi^2(d/2)^2/\eta} = e^{-m^2/2}$. For the remainder Gaussians to
be truncated at the same level, we need that
\begin{align}
  e^{-\xi^2(\deltaL/2)^2/(1-\eta)} \le e^{-m^2/2},
\end{align}
i.e.
\begin{align}
  \deltaL \ge \sqrt{2 (1-\eta) m^2/\xi^2} .
\end{align}
To guarantee that both Gaussians have proper support, we thus need
\begin{align}
  \deltaL \ge 
  \begin{cases}
    d & \text{ if } \eta \ge 1 ,\\
    \max\left(d, \sqrt{2 (1-\eta) m^2/\xi^2} \right) & \text{ if } \eta < 1 .
  \end{cases}
  \label{eq:deltaL_cond}
\end{align}
With this extra support for $\eta < 1$ the approximation errors are
decoupled from the Fourier space truncation errors, which are further
discussed in section \ref{sec:fourier-space-error}. An example of this
decoupling is shown in Figure \ref{fig:P_errors}, where it can be seen
that the larger choice of $\deltaL$ is actually only needed if the grid size $\Mtilde$ is picked
larger than necessary for a given error tolerance.

\section{Evaluating the real space component}
\label{sec:eval-real-space} 
The real space part of the free-space Ewald sum \eqref{eq:usumEwald}
has the general form
\begin{align}
  \v u^R(\xb) = \sum_{\ni=1}^N \Green^R(\xb - \xb_\ni) \cdot
  \fb(\xb_\ni).
  \label{eq:usum_r}
\end{align}
Since $\Green^R(r)$ decays rapidly (roughly as $e^{-\xi^2 r^2}$), the
sum can be truncated outside some truncation radius $r_c$. Assume that
we wish to evaluate the potential at points $\xb_\mi$, $\mi=1\dots
N$. The expression that we need to evaluate is then
\begin{align}
  \v u^R(\xb_\mi) = \sum_{\substack{\ni=1\\|\xb_\mi - \xb_\ni| \le r_c}}^N
  \Green^R(\xb_\mi - \xb_\ni) \cdot \fb(\xb_\ni), \quad \mi=1\dots N.
\end{align}
Naively implemented, this has an $\ordo(N^2)$ computational cost. It
is however straightforward to find the interaction list of each target
point $\xb_\mi$ by first creating a cell list \cite{Frenkel2001}. This
reduces the real space cost to $\ordo(N)$, under the assumption that
the average number of interactions of each target point stays
constant when $N$ changes.

\section{Truncation errors}
\label{sec:truncation-errors}

Truncation errors are introduced when we cut off the real space
interactions outside a radius $r_c$, and when we truncate the Fourier
space integral outside a maximum wave number $\kmax$. The magnitudes
of these errors can be accurately estimated through the analysis
methodology introduced by Kolafa \& Perram \cite{Kolafa1992} for
periodic electrostatic force computations. Denoting by $\Delta \v u(\v
x)$ the pointwise truncation error in the solution, one can then
derive statistical error estimates for the root mean square (RMS)
truncation error, defined as
\begin{align}
  \delta \v u = \sqrt{ \frac{1}{N} \sum_{\ni=1}^N \left|\Delta \v u(\xb_{\ni})\right|^2 }.
\end{align}
The analyses for both the real and Fourier space components rely on
the following property:
\begin{lemma}{(Kolafa \& Perram \cite[appx.A]{Kolafa1992})}
  Let $(\xb_{\ni}, q_n)$ be a configuration of point sources, and let
  \begin{align}
    \Delta F(\xb) = \sum_{\ni=1}^N q_n \Delta f(\xb-\xb_{\ni}),
  \end{align}
  be an error measure due to a set of pointwise errors, $\Delta f =
  f_{\text{approx}}-f_{\text{exact}}$. Assuming that the points are
  randomly distributed, and that $\Delta F$ has a Gaussian
  distribution, then the root mean square (RMS) error
  \begin{align}
    \delta F = \sqrt{\frac{1}{N} \sum_{\ni=1}^N \left(\Delta F(\xb_{\ni})\right)^2}
  \end{align}
  can be approximated as
  \begin{align}
    \delta F^2 \approx  \frac{1}{|V|} \sum_i q_i^2 \int_V (\Delta f(\v r))^2 \dif \v r,
  \end{align}
  where $V$ is the volume enclosing all point-to-point vectors $\v
  r_{ij} = \xb_i - \xb_j$.
  \label{lemma:rms}
\end{lemma}

\subsection{Fourier space truncation error}
\label{sec:fourier-space-error}

The Fourier space error comes from truncating the integral of the
Fourier transform outside a maximum wave number $\kmax$,
\begin{align}
  \Delta\v u^F(\xb) = \frac{1}{(2\pi)^3} \int_{k > \kmax} 
  \widehat G^F(\v k, \xi) \cdot \sum_{\ni=1}^N \v f(\xb_{\ni}) 
  e^{i\v k\cdot(\xb - \xb_{\ni})} \dif \v k, 
\end{align}
where all points $\xb_{\ni}$ are contained in a cube of size $L$. 
In our case the integral is approximated using an FFT over an $M^3$
grid covering an $L^3$ domain, such that
\begin{align}                   
  \kmax = \frac{2\pi}{\Ltilde} \frac{\Mtilde}{2} .
\end{align}
The RMS of the truncation error is given by
\begin{align}
  \delta\v u^F =\sqrt{\frac{1}{N} \sum_{\ni=1}^N \left| \Delta\v u^F(\xb_{\ni}) \right|^2 },
\end{align}
and can be estimated using the method of Kolafa \& Perram. Such
estimates already exist for the periodic stokeslet \cite{Lindbo2010}
and rotlet \cite{AfKlinteberg2016rot} potentials, as well as for the
Beenakker decomposition of the stresslet
\cite{AfKlinteberg2014a}. However, it turns out that periodic
estimates fail for free-space potentials that are based on the
truncated biharmonic potential $\biharmonic^\trunc$. This is because
the dominating term of $\widehat\biharmonic^\trunc$ for large $k$ is a
factor $\trunc^2 k^2/2$ larger than $\widehat\biharmonic$. For
potentials based on the truncated harmonic potential
$\harmonic^\trunc$ one can use the periodic estimates, as the
difference in magnitude between $\widehat\harmonic^\trunc$ and
$\widehat\harmonic$ is negligible. We can thus use existing estimates
for the rotlet, while we need to derive new ones for the stokeslet and
stresslet. The final set of estimates is shown in Table \ref{tab:est_fourier}.

\subsubsection{Stokeslet}
Beginning with the stokeslet potential, we consider the truncation
error contribution from a single source located at the origin,
\begin{align}
  \Delta u^F_j(\xb) = e_{jl}(\xb)f_l,
\end{align}
where
\begin{align}
  e_{jl}(\v r) &= \frac{1}{(2\pi)^3} \int_{k > \kmax} \left(1 +
    \frac{k^2}{4\xi^2}\right) \widehat\biharmonic^\trunc(\v k) e^{-(k/
    2\xi)^2} e^{i \v k\cdot \v r} k^2 \left(\delta_{jl} - \hat k_j
    \hat k_l \right) \dif \v k,
\end{align}
and
\begin{align}
  \widehat \biharmonic^\trunc(\v k) = 4 \pi \frac{
    (2-\trunc^2k^2)\cos(\trunc k) + 2 \trunc k \sin(R k) - 2
  }{k^4} .
\end{align}
We now keep only the highest order term in $k$, which dominates the
error for large $\kmax$,
\begin{align}
  e_{jl} &\approx -\frac{4 \pi \trunc^2}{4\xi^2(2\pi)^3} \int_{k >
    \kmax} k^2\cos(\trunc k) e^{-(k/ 2\xi)^2} e^{i \v k\cdot \v r}
  \left(\delta_{jl} - \hat k_j \hat k_l \right) \dif \v k.
\end{align}
The error should be independent of the coordinate system orientation,
and since we are deriving a statistical error measure we approximate
the directional component by its root mean square (computed using
spherical coordinates),
\begin{align}
  \left(\delta_{jl} - \hat k_j \hat k_l \right) \approx
  \sqrt{\frac{1}{9}\sum_{j,l=1}^3 \left(\delta_{jl} - \hat k_j \hat
      k_l \right)^2} = \frac{\sqrt{2}}{3}.
\end{align}
We now integrate in spherical coordinates, choosing a system such that
$\hat k_3 \parallel \hat z$. Integration in $\theta\in[0,\pi]$ then gives us
\begin{align}
  e_{jl} &\approx 
-\frac{\sqrt{2}}{3} \frac{\trunc^2}{2\xi^2\pi r}
  \int_{k > \kmax} k^3 \cos(\trunc k) \sin(rk)
  e^{-(k/ 2\xi)^2}
  \dif k.
\end{align}
To get a good approximation of $e_{jl}$, we need to approximate the
remaining integral. 
First, we have from the exponential decay that
the dominating contribution will come from the beginning of the
interval, where $k \approx \kmax$. This allows us to approximate
\begin{align}
  I &= \int_{k > \kmax} k^3 \cos(\trunc k) \sin(rk) e^{-(k/ 2\xi)^2}
  \dif k \approx \kmax^3 \int_{k > \kmax}  \cos(\trunc k) \sin(rk) e^{-(k/ 2\xi)^2} \dif k .
\end{align}
Next, we have that $\trunc \gg r$, so $\trunc$ is the dominating
frequency in the integrand, such that we can write
 $\cos(\trunc k)\sin(rk) \approx \cos(\trunc k) \sin(r\kmax)$, 
and
\begin{align}
  I \approx \kmax^3 \sin(\kmax r) \int_{k > \kmax}  e^{i\trunc k-(k/ 2\xi)^2} \dif k,
\end{align}
where we implicitly assume that the real part of the complex
exponential is our quantity of interest. A final approximation (again
assuming $k\approx \kmax$) makes this integrable, and we get
\begin{align}
|  I | &\approx \left| \frac{\kmax^3 \sin(\kmax r)}{i\trunc - \kmax/ 2\xi^2}
  \int_{k > \kmax} (i\trunc - k/ 2\xi^2)
  e^{iRk-(k/ 2\xi)^2} \dif k \right| \\
  &= \left| \frac{\kmax^3 \sin(\kmax r)}{i\trunc - K/ 2\xi^2}
    e^{i\trunc\kmax-(\kmax/ 2\xi)^2} \right| 
\le \frac{\kmax^3 |\sin(\kmax r)|}{|i\trunc - \kmax/ 2\xi^2|}
  e^{-(\kmax/ 2\xi)^2}.
\end{align}

We have (for a cube) that $\trunc \ge \sqrt{3} \Mtilde h$, while
$\kmax = \pi \Mtilde/\Ltilde = \pi / h$. Typical parameter values are
around $\kmax/\xi = \mathcal O(10)$ and $\Mtilde = \ordo(50)$, so
$\kmax/2\xi^2 = \mathcal O(50 h/\pi) \ll \trunc$ and $|i\trunc -
\kmax/ 2\xi^2| \approx \trunc$. We can therefore write
\begin{align}
  | I | \approx \frac{\kmax^3 |\sin(\kmax r)|}{\trunc} e^{-(\kmax/ 2\xi)^2},
\end{align}
and
\begin{align}
  |e_{jl}| \approx \frac{\sqrt{2}}{6} \frac{\trunc \kmax^3 |\sin(\kmax
    r)|}{\xi^2\pi r} e^{-(\kmax/ 2\xi)^2} .
\end{align}
We can now use Lemma \ref{lemma:rms} to estimate the statistical error
by integrating over a sphere of radius $L/2$ (which then contains all
point sources),
\begin{align}
  \left( \delta\v u^F \right)^2 
  &\approx \sum_{\ni=1}^N \sum_{j=1}^3 f_l^2(\xb_{\ni}) \frac{1}{|V|} \int_V e_{jl}^2(r) \dif \v r \\
  &\approx Q\frac{6}{\pi L^3} 3 \left( \frac{\sqrt{2}}{6} \frac{\trunc
      \kmax^3}{\xi^2\pi} e^{-(\kmax/ 2\xi)^2} \right)^2 \int_0^{L/2}
  \sin^2(\kmax r) 4\pi \dif r,
\end{align}
where
\begin{align}
  Q = \sum_{\ni=1}^N |\v f(\xb_{\ni})|^2 .
  \label{eq:Q_stokeslet}
\end{align}
Assuming that $\sin^2(\kmax r)$ has many oscillations in the interval
$[0,L/2]$, we replace it by its average value $1/2$, such that
 $\int_0^{L/2} \sin^2(\kmax r) 4\pi \dif r \approx \pi L$.
Finally, we can write the stokeslet truncation error estimate as
\begin{align}
  \delta\v u^F \approx \sqrt{Q}\frac{\trunc \kmax^3}{\xi^2 \pi L} e^{-(\kmax/
    2\xi)^2} .
\end{align}

\subsubsection{Stresslet} 
For the stresslet, the derivation for the error estimate is completely
analogous to the one for the stokeslet. The difference is that the
leading order term is $k^3$ instead of $k^2$, and that the RMS of the
directional component is
\begin{align}
  \sqrt{
  \frac{1}{27} \sum_{j,l,m=1}^3 \left(
    \left(\delta_{jl}\hat k_m+\delta_{lm}\hat
      k_j+\delta_{mj}\hat k_l\right) - 2\hat k_j\hat k_l\hat k_m
  \right)^2}
  = \frac{7}{27} .
\end{align}
This allows us to directly write the stresslet truncation error
estimate as
\begin{align}
  \delta\v u^F \approx \sqrt{\frac{7Q}{6}}\frac{\trunc \kmax^4}{\xi^2 \pi L} e^{-(\kmax/
    2\xi)^2},
\end{align}
where
\begin{align}
  Q = \sum_{\ni=1}^N \sum_{l,m=1}^3 q_l^2(\xb_{\ni}) n_m^2(\xb_{\ni}).
  \label{eq:Q_stresslet}
\end{align}

\begin{table}[htbp]
  \centering
  \begin{tabular}{c|c|c|c}
    & Stokeslet, $\stokeslet^F$ 
    & Stresslet, $\stresslet^F$ 
    & Rotlet, $\rotlet^F$ 
    \\ \hline
    &&&\\
    $\delta\v u^F$ & 
    $\displaystyle \sqrt{Q}\frac{\trunc \kmax^3}{\xi^2 \pi L} e^{-\kmax^2/4\xi^2}$ &
    $\displaystyle \sqrt{\frac{7Q}{6}}\frac{\trunc \kmax^4}{\xi^2 \pi L} e^{-\kmax^2/4\xi^2}$ &
    $\displaystyle \sqrt{\frac{8\xi^2Q}{3\pi L^3 \kmax}} e^{-\kmax^2/4\xi^2}$ 
    \\
    &&&\\
  \end{tabular}
  \caption{Fourier space truncation errors for the stokeslet, stresslet, and rotlet \cite{AfKlinteberg2016rot}. The quantity $Q$ is defined as in \eqref{eq:Q_stokeslet} for the stokeslet and rotlet, and as in \eqref{eq:Q_stresslet} for the stresslet.}
  \label{tab:est_fourier}
\end{table}

\begin{figure}[htbp]
  \centering
  \begin{subfigure}[b]{0.32\textwidth}
    \centering
    \includegraphics[width=\textwidth]{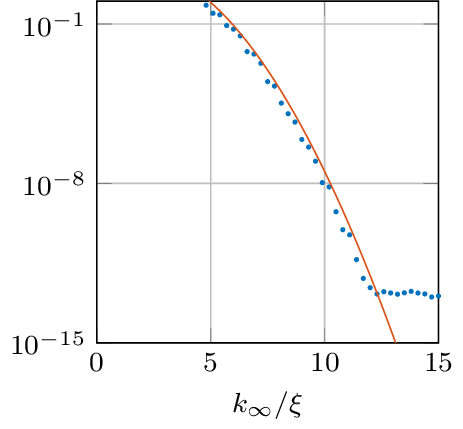}
    \caption{Stokeslet}
    \label{fig:fourier_err:stokeslet}
  \end{subfigure}  
  \begin{subfigure}[b]{0.32\textwidth}
    \centering
    \includegraphics[width=\textwidth]{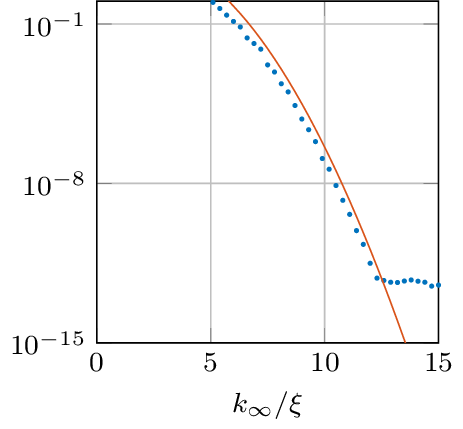}
    \caption{Stresslet}
    \label{fig:fourier_err:stresslet}
  \end{subfigure}  
  \begin{subfigure}[b]{0.32\textwidth}
    \centering
    \includegraphics[width=\textwidth]{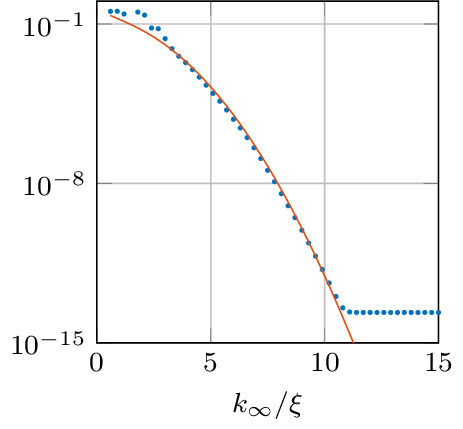}
    \caption{Rotlet}
    \label{fig:fourier_err:rotlet}
  \end{subfigure}  
  \caption{RMS of relative Fourier space truncation errors for the
    stokeslet, stresslet and rotlet. Dots are measured value, solid
    lines are computed using the estimates of Table
    \ref{tab:est_fourier}. The system is $N=10^4$ randomly distributed
    point sources in a cube with sides $L=3$, with $\kmax=\pi \Mtilde
    /\Ltilde$, $\xi=3.49$, $M=1\dots50$, and $P=32$.}
  \label{fig:fourier_err}
\end{figure}

\subsection{Real space truncation error}

The real space truncation error is due to neglecting interactions in
the real space sum for which $r > r_c$, 
\begin{align}
  \Delta\v u^R(\xb) = \sum_{|\xb - \xb_{\ni}| > r_c} %
  G^R(\xb - \xb_{\ni}) \cdot \v f(\xb_{\ni}) .
\end{align}
The RMS of this error is given by
\begin{align}
  \delta\v u^R =\sqrt{\frac{1}{N} \sum_{\ni=1}^N \left| \Delta\v u^R(\xb_{\ni}) \right|^2 }
\end{align}
Following the analysis by Kolafa \& Perram, we can use Lemma
\ref{lemma:rms} to estimate $\delta\v u^R$ as
\begin{align}
  \left( \delta\v u^R \right)^2 \approx \frac{1}{L^3}\sum_{\ni=1}^N 
  \left( \v f(\xb_{\ni}) \right)^2
  \cdot \int_{r>r_c} \left(G^R(\v r)\right)^2 \dif \v r .
\end{align}
Estimates based on this approximation are already available in the
literature for the stokeslet \cite{Lindbo2011e} and rotlet
\cite{AfKlinteberg2016rot} decompositions used in this paper, and are
shown in the summary in Table \ref{tab:est_real}. We will here derive
a similar estimate for the Hasimoto decomposition of the stresslet,
essentially by repeating the derivation of \cite{AfKlinteberg2014a}
for the Beenakker decomposition.

The real space component of the stresslet has the form
\begin{align}
  \stresslet^R_{jlm}(\xi, \v r) = A_1(\xi, r) \hat r_j \hat r_l \hat r_m +
  A_1(\xi, r) (\delta_{jl}\hat r_m + \delta_{lm}\hat r_j + \delta_{mj}\hat r_l ),
\end{align}
and the RMS error is approximated as
\begin{align}
  \left( \delta\v u^R \right)^2 \approx \frac{1}{L^3}\sum_{\ni=1}^N
  \sum_{j=1}^3
  q_l^2(\xb_{\ni}) q_m^2(\xb_{\ni}) \int_{r>r_c} \left(T^R_{jlm}(\v
    r)\right)^2 \dif \v r .
  \label{eq:TRrms}
\end{align}
Arguing that the error should be independent of the coordinate system
orientation, we replace $(\stresslet^R_{jlm})^2$ by its average value
over the tensor components, computed using spherical coordinates
\begin{align}
  \sum_{j=1}^3 \left(T^R_{jlm}(\v r)\right)^2 \approx%
  3 \overline{\left(\stresslet^R\right)^2} = \frac{3}{27}
  \sum_{j,l,m=1}^3 \left( T_{jlm}^R \right)^2 = \frac{1}{9}
  \left(A_1^2+6 A_1 A_2+15 A_2^2\right) .
  \label{eq:TR2avg}
\end{align}
This quantity has only radial dependence and is integrable, 
\begin{align}
  \begin{split}   
    \int_{r > r_c} 3 \overline{\left(\stresslet^R\right)^2} 4\pi r^2 \dif r =& 
    -\frac{32}{3} \sqrt{\pi } \xi e^{-\xi ^2 r_c^2} \text{erfc}\left(\xi
      r_c\right)+21 \sqrt{2 \pi } \xi \text{erfc}\left(\sqrt{2} \xi
      r_c\right)\\&+\frac{16 \pi \text{erfc}\left(\xi
        r_c\right){}^2}{r_c}+\frac{4}{9} \xi ^2 r_c e^{-2 \xi ^2 r_c^2}
    \left(28 \xi ^2 r_c^2-3\right)
    \\ \approx & 
    \frac{112}{9} \xi ^4 r_c^3 e^{-2 \xi ^2 r_c^2},
  \end{split}
  \label{eq:TR2int}
\end{align}
where we have kept only the dominating term (for large $\xi r_c$) in
the last step. Combining \eqref{eq:TRrms}, \eqref{eq:TR2avg} and
\eqref{eq:TR2int} gives us the estimate for the stresslet shown in
Table \ref{tab:est_real}.

\begin{table}[htbp]
  \centering
  \begin{tabular}{c|c|c|c}
    & Stokeslet, $\stokeslet^R$ 
    & Stresslet, $\stresslet^R$ 
    & Rotlet, $\rotlet^R$ 
\\
    \hline
    &&&\\
    $\delta\v u^R$ 
    & $\displaystyle \sqrt{\frac{4 Q r_c}{L^3}} e^{-\xi^2 r_c^2}$
    & $\displaystyle \sqrt{\frac{112 Q \xi^4 r_c^3}{ 9 L^3 }} e^{-\xi^2 r_c^2}$ 
    & $\displaystyle \sqrt{\frac{8 Q}{3 L^3 r_c}} e^{-\xi^2 r_c^2}$ 
    \\
    &&&\\
  \end{tabular}
  \caption{Real space truncation errors for the stokeslet \cite{Lindbo2011e}, stresslet, and rotlet \cite{AfKlinteberg2016rot}. The quantity $Q$ is defined in the same way as for the Fourier component.}
  \label{tab:est_real}
\end{table}

\begin{figure}[htbp]
  \centering
  \begin{subfigure}[b]{0.32\textwidth}
    \centering
    \includegraphics[width=\textwidth]{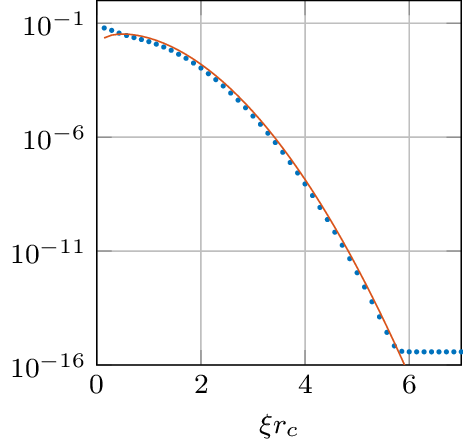}    
    \caption{Stokeslet}
    \label{fig:real_err:stokeslet}
  \end{subfigure}  
  \begin{subfigure}[b]{0.32\textwidth}
    \centering
    \includegraphics[width=\textwidth]{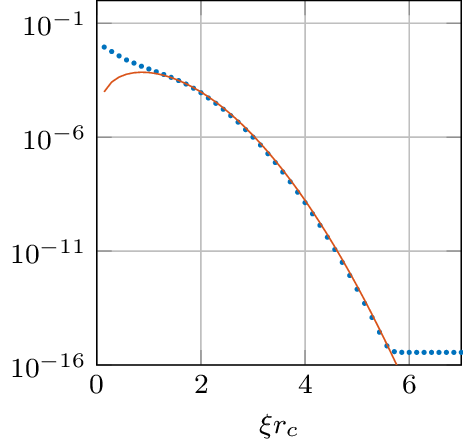}
    \caption{Stresslet}
    \label{fig:real_err:stresslet}
  \end{subfigure}  
  \begin{subfigure}[b]{0.32\textwidth}
    \centering
    \includegraphics[width=\textwidth]{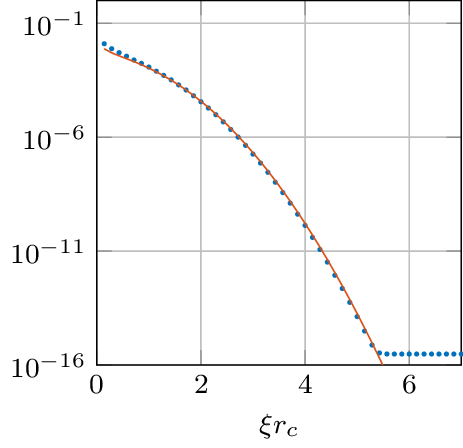}    
    \caption{Rotlet}
    \label{fig:real_err:rotlet}
  \end{subfigure}  
  \caption{RMS of relative real space truncation errors for the
    stokeslet, stresslet and rotlet. Dots are measured value, solid
    lines are computed using the estimates of Table
    \ref{tab:est_real}. The system is $N=2000$ randomly distributed
    point sources in a cube with sides $L=3$, with $\xi=4.67$, and
    $r_c \in [0, L/2]$.}
  \label{fig:real_err}
\end{figure}

\section{Summary of method}
\label{sect:summary}

We now summarize the free-space fast Ewald method for Stokes
potentials. Our goal is to compute a discrete-sum potential of the
type \eqref{eq:usum},
\begin{align}
  \v u(\xb) = \sum_{\ni=1}^N \Green(\xb - \xb_{\ni}) \cdot
  \fb(\xb_{\ni}),
\end{align}
for a set of $N$ target point $\xb$, with $\Green$ being the stokeslet
\eqref{eq:def_stokeslet}, stresslet \eqref{eq:def_stresslet} or rotlet
\eqref{eq:def_rotlet}. We assume that all target and source points are
contained in the cubic domain $\domain = [0,L]^3$.

Using an Ewald decomposition (sec. \ref{sec:ewald-free-space}) and
an Ewald parameter $\xi > 0$, we split the potential into a
short-range part $\v u^R$ acting locally, and a long-range part $\v
u^F$ computed in Fourier space,
\begin{align}
  \v u(\xb) = \v u^R(\xb, \xi)+ \v u^F(\xb, \xi) + \v u^{\text{self}}(\xb, \xi),
\end{align}
where $\v u^{\text{self}}$ refers to the self-interaction term
\eqref{eq:Slimit}, which has to be taken into account only for the
stokeslet potential.

The real space component is truncated outside an interaction radius
$r_c$,
\begin{align}
  \v u^R(\xb) \approx \sum_{\substack{\ni=1\\|\xb - \xb_\ni| \le r_c}}^N
  \Green^R(\xb - \xb_\ni) \cdot \fb(\xb_\ni).
\end{align}
This is a local operation in the neighborhood of each target point
$\xb$, and can be efficiently evaluated using e.g. a cell list
(sec. \ref{sec:eval-real-space}).

The Fourier space component is evaluated through a Fourier integral,
truncated at a maximum wave number $\kmax$,
\begin{align}
  \v u^F(\xb,\xi) \approx \frac{1}{(2\pi)^3} \int_{|\kb| \le \kmax} e^{i\v
    k\cdot\xb} 
  \widehat{\Green}^{F,\trunc}(\v k, \xi)
  \cdot \sum_{\ni=1}^N \fb(\xb_{\ni}) e^{-i\v
    k \cdot \xb_{\ni}} \dif \v k.
\end{align}
The superscript $\trunc$ denotes that we have removed the singularity
in the integrand (at $k=0$) by truncating the original Green's
function outside a maximum interaction radius $\trunc$
(eq. (\ref{eq:AS}--\ref{eq:AR}), eq. \eqref{eq:Gwoexp}, sec.
\ref{sec:free-space-solution}). The integral is evaluated using the
spectral Ewald method (sec. \ref{sec:discretization}), which uses FFTs
on an $\Mtilde^3$ grid to efficiently compute the long-range
interactions. The method requires the domain $\domain$ to be extended
by a length $\deltaL$ \eqref{eq:deltaL_cond} for accurate function
support, and then zero-padded by a factor 2 for the convolution to be
aperiodic when using FFTs. In fact, an oversampling factor $\osamp \ge
1+\sqrt{3} \approx 2.8$ is required to accurately resolve the
truncated Green's function (sec. \ref{sec:oversampling}), but the cost
of that can be reduced to a precomputation step involving the
fundamental solution to the harmonic or biharmonic equation
(sec. \ref{sec:precomputation}).

\subsection{Computational complexity}
\label{sec:comp-complexity}

We wish to express how the computational cost of the method scales
with an increased number of sources and targets $N$, which we assume
to be evenly distributed in the domain $\domain$. The system can be
scaled up in two different ways: by increasing the point density in a
fixed domain, or by increasing the domain size $L$ with a fixed point
density. Either way, the scaling arguments have as their starting
point that the real space sum be $\ordo(N)$. This is achieved by
keeping a constant number of near neighbors (within $r_c$) for each
target under scaling. Additionally, we want the level of the
truncation errors to be constant, which is achieved by keeping $\xi
r_c$ and $\Mtilde \xi^{-1} L^{-1}$ constant.

If $N$ increases with $\domain$ fixed, then $r_c \propto N^{-1/3}$ is
required for an $\ordo(N)$ real space sum. If the accuracy is to
remain constant, then $\xi \propto r_c^{-1} \propto N^{1/3}$ and the
grid size is scaled as $\Mtilde \propto \xi \propto N^{1/3}$. This puts
the Fourier space cost at $\ordo(\Mtilde^3 \log \Mtilde) \propto
\ordo(N \log N)$.

If the domain size $L$ increases with a fixed point density, then $N
\propto L^{3}$ and the real space sum is $\ordo(N)$ if we keep $r_c$
and $\xi$ constant. Then $\Mtilde \propto L \propto N^{1/3}$, such
that the Fourier space cost is $\ordo(\Mtilde^3 \log \Mtilde) \propto
\ordo(N \log N)$.

\subsection{Parameter selection}
\label{sec:parameter-choice}

For a given system ($N$ charges in a domain of size $L$), the required
parameters for our free-space Ewald method are the Ewald parameter
$\xi$, the real space truncation radius $r_c$, the number of grid
points $M$ covering the original domain, and the Gaussian support
width $P$. Based on these parameters one can then set $\deltaL$ using
\eqref{eq:deltaL_cond}, which then gives $\Ltilde = L + \deltaL$. This
in turn gives $\Mtilde$, by satisfying $h=L/M=\Ltilde/\Mtilde$. We
will here draft a strategy for optimizing $\xi$, $r_c$, $M$ and $P$ in
a large-scale numerical computation.

For a given value of $\xi$ and absolute error tolerance $\epsilon$,
close-to optimal values for $M$ and $r_c$ can be computed using the
estimates in Tables \ref{tab:est_fourier} and \ref{tab:est_real}. The
support width $P$ affects the error in the Fourier space component,
and we have in practice observed that for the relative error, $P=16$
gives 8 digits of accuracy, $P=24$ gives 12 digits, and $P=32$ is
enough to guarantee that the approximation errors are at roundoff.

Which value of $\xi$ to choose is highly implementation dependent, as
the variable is used to shift the workload between the real and
Fourier space components. A straightforward strategy for finding an
optimal value is to start with a small but representative subset of
the original system, and compute a reference solution for that
subset. Picking a starting value for $\xi$, one then sets $P=32$, and
adjusts $M$ and $r_c$ until the error tolerance is strictly met. Then
$P$ can be decreased in steps of two\footnote{$P$ a multiple of two is
  favorable for code optimization using vector instructions.} until
the tolerance is reached again. Using this starting point for $(\xi,
r_c, M, P)$, one then does a parameter sweep in $\xi$ for finding the
configuration with the smallest runtime, while keeping $\xi r_c$ and
$M/\xi$ constant during the sweep. Once an optimal setup is found, the
original (large) system can be computed using the same set of
parameters, except $M$ which is scaled such that $L/M$ remains
constant for both systems.

\section{Results}
\label{sect:numres}

We consider systems of $N$ random point sources drawn from a
uniform distribution in a box of size $L^3$. 
We evaluate the sum \eqref{eq:usum} with \textit{stokeslets}
\eqref{eq:def_stokeslet}, \textit{stresslets} \eqref{eq:def_stresslet}
and \textit{rotlets} \eqref{eq:def_rotlet} using our free space
Spectral Ewald (FSE) method, at the same $N$ target locations.  
All components of the force/source strengths are random numbers from a
uniform distribution on $[-1,1]$. 
All computationally intensive routines are written in C and are called
from Matlab using MEX interfaces.  The results are obtained on a
desktop workstation with an Intel CoreTM i7-3770 Processor (3.40 GHz)
with four cores and 8 GB of main memory.
To measure the actual errors, we compare to the result from evaluating
the sum by direct summation. 

\begin{figure}[htbp]
  \centering
  \includegraphics[width=0.49\textwidth,height=0.5\textwidth]{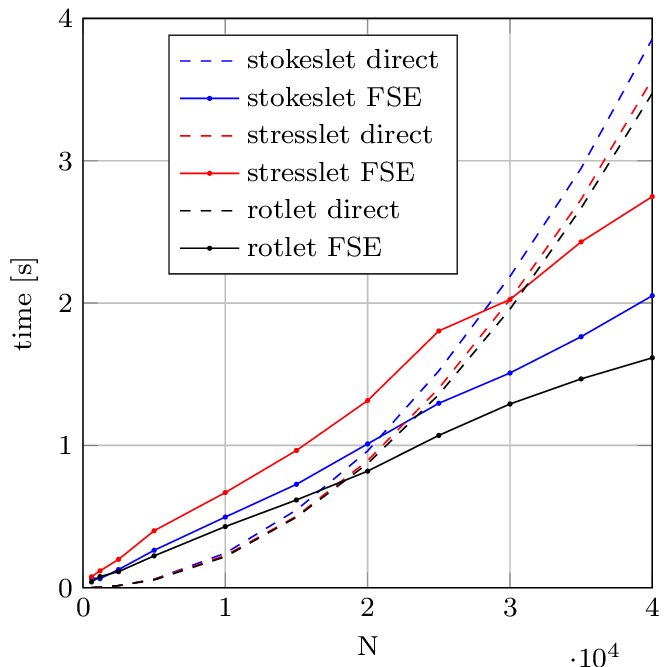}  
  \includegraphics[width=0.49\textwidth,height=0.5\textwidth]{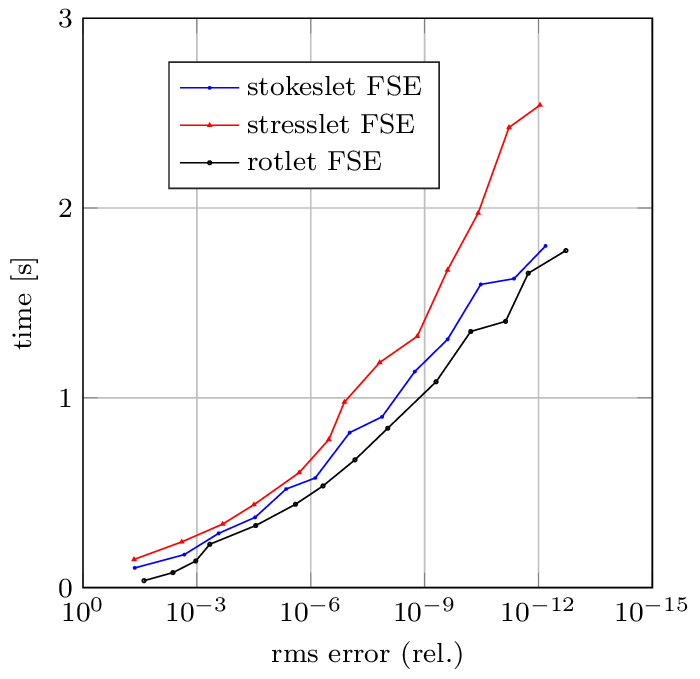} 
  \caption{(left) Comparison of direct and fast evaluation of the sum in
    \eqref{eq:usum} for the Stokeslet, Stresslet, and Rotlet 
    including the precomputation step. 
The system grows at constant density ($N/L^3$ is constant), with
$\xi=7$ and $r_c$ constant for all values of $N$. 
The relative  rms error is less than $0.5\times10^{-8}$.
(right) Runtime of computing \eqref{eq:usum} using FSE for  all kernels as a function of the relative rms error. $N=20 \, 000$, $L=2$ and $\xi=7$.
}
  \label{fig:runtime_all_kernels}
\end{figure}

In the left plot of figure \ref{fig:runtime_all_kernels}, the computing time
for evaluation of the sums is plotted versus $N$, for all three
kernels and for both the Spectral Ewald (FSE) method and direct
summation.  The parameters in the Spectral Ewald method have been set
to keep the relative rms error below $0.5 \cdot 10^{-8}$.
The optimal value of $\xi$ cannot be determined theoretically, since
it is implementation and hardware dependent. When we vary $N$ in
figure \ref{fig:runtime_all_kernels}, we change the size of the box,
to keep a constant number density $N/L^3=2500$. If an optimal value of
$\xi$ is determined for one system (see discussion in section
\ref{sec:parameter-choice}), the same value can be kept as the system
is scaled up or down in this manner.
The parameters $r_c$, $P$ and the grid resolution $L/M$ are kept
constant as $N$ and hence $L$ is increased, yielding an increase in
the grid size. 
We have used $\xi=7$ for all three kernels, $r_c=0.63$, $0.63$ and
$0.58$ for the stokeslet, stresslet and rotlet, respectively, and
$P=16$ for all kernels. For $L=2$, $M$ is set to $48$, $50$ and $38$
for the three kernels, and is then scaled with $L$. 

The precomputation step does not depend on the location
of the sources, and can be performed once the size of the domain is
set. 
The precomputation cost can therefore
usually be amortized over many calls to the method, as a simulation
code is run for many time steps and possibly iterations within time steps. 
Despite this, we have chosen to plot the runtimes including the
precomputation cost, and later discuss it in more detail. 
From these data (figure \ref{fig:runtime_all_kernels} left), we can
find the approximate breakeven points, i.e. the values of $N$ for
which any larger system will benefit from using the fast method.  We
find it to be approximately $N=20 \, 000$ for the stokeslet,
$30 \, 000$ for the stresslet and $19 \, 000$ for the rotlet with
precomputation, which is reduced to $18 \, 000$, $25 \, 000$ and
$15 \, 000$ without the precomputation step.  If the precomputation
step is to be done only once, the decomposition parameter $\xi$ should
however be chosen differently for optimal performance, which would
bring down the break even point further.  Note that this is a strict
error tolerance. For lower accuracy requirements, the cross over
occurs at lower values of $N$.  These are higher values than have
previously been reported in the literature, e.g. in
\cite{Tornberg2008}, where $N=5000$ was reported as the breakeven for
the Stokelet. There are two factors affecting these numbers, one is
that these results are run on multiple cores for which the direct sum
parallelizes better than the FFTs involved in the fast method. The
other factor is that the direct sums relatively speaking have become
faster to evaluate also on a single core, where compilers can speed up
the code significantly using vector instructions, while the more
complicated algorithms cannot benefit from this as extensively.

To make sure that our method is competitive, we have compared to a
fast multipole implementation available as free software
\cite{fmm_webpage}, running both codes on a single core and comparing
timings for the Stokeslet.  We set the accuracy level to six digits
in the FMM. For $N=20 \, 000$, this yields a relative rms error of
about $5.6 \cdot 10^{-9}$, and we set the parameters for the FSE
method to obtain a similar error level (for this case we get
$4.3 \cdot 10^{-9}$).  For $N=20 \, 000$, the FSE code (including
precomputation) and the FFM code both use about 3 seconds.  The direct
evaluation takes 3.6s with our code and 6.7s with the code provided
with the FMM package.  It should however be noted, that the FMM as
well as the direct code from that package, returns not only the three
vector components produced by the Stokeslet, but also the associated
(scalar) pressure, which increases the cost somewhat.  The breakeven
point for both FSE and FMM is about $N=17 \, 000$ when comparing to
our direct code. If we instead compare to the direct code in the FMM
package, the break even point for the FMM decreases to $N=10 \, 000$.
Most fair would be to compare the FMM to a direct sum written as the
faster one, but including also the pressure component, which should
place the break even point between the two numbers above.  For the FSE
code, assuming that the precomputation will be done only once, and
choosing $\xi$ instead to optimize the runtime without precomputation,
the break even point drops from $N=17 \, 000$ to $11 \, 000$.

Let us consider also a larger system with $N=400 \, 000$, with
$L=5.4288...$ such that $N/L^3=2500$. For the stokeslet summation by
FSE, we pick the parameters $\xi=8$, $r_c=0.5651$, $M=144$ and $P=16$
to obtain a relative rms error of $5 \cdot 10^{-8}$. This means that
the FFTs are computed for grids with $\Mtilde=2(M+P)$.  The time for
evaluation is about 64 seconds (including the precomputation), and the
speed-up compared to our direct evaluation of the sum is a factor of
about $23$.  Excluding the precomputation cost, the computing time is
reduced by 15 s, and this factor increases to $29$.  For the FMM, the
evaluation time is about 180 seconds, yielding a speed-up of a factor
of about $8$ compared to our direct sum or a factor of $15$ as
compared to the one provided with the FMM code \cite{fmm_webpage}.
Checking the relative rms error from both the FSE and FMM
computations, they are similar, around $0.5 \cdot 10^{-8}$ for FSE and
$10^{-8}$ for the FMM.  Hence, for this example on a single core, the
FSE method including the precomputation is almost three times as fast
as the FMM method, but the difference would be reduced somewhat if the
time for computing the extra pressure component was excluded.

In the adaptive FMM code, a box is split into 8 children boxes if the
number of sources is larger than a set value. If any of the children
boxes still have too many source points, it is split again. With a
uniform distribution of points, most leaf boxes are on the same level
of refinement, which in this case will be four divisions.  The curve
for computational cost versus N will not be smooth, since this is a
discrete process (either you keep the box as one, or you split into
eight), which changes the cost balance between different parts of the
algorithm. This is why the larger computational cost for the FMM
method in this case could not be predicted considering the timing for
$N=20 \, 000$, where the timing of the FMM and FSE methods were
similar.

We did not set out to make a thorough comparison of the two
methods. All results are for uniform distributions of source points.
Typically the FSE performs better compared to the FMM for higher
accuracies. Moving towards an increasingly non-uniform point
distribution, the adaptivity of the FMM will at some point pay off.
With this, we have however showed that the FSE method is competitive with
the FMM.

To show how the computational cost depends on the accuracy
requirements, we now consider a fixed system with $N=20 \, 000$ sources in
a box with $L=2$ and vary the error tolerance.
In the right plot of figure  \ref{fig:runtime_all_kernels}, we plot the runtime for
summing the stokeslet, stresslet and rotlet kernels as a function of
the relative rms error in the result.  Computing the $k$-space
contribution for the stokeslet and rotlet involves gridding of three
vector components, three FFTs, a scaling in Fourier-space, three
inverse FFTs and the quadrature step for the three components of the
solution, see the algorithm in section \ref{sec:discretization}. The
stresslet instead requires the gridding of 9 components and hence 9
FFTs. After the scaling step, there are three resulting vector
components, as for the other kernels.  All three kernels require the
same amount of precomputing.  Hence, it is not surprising to see that
the stresslet is the most expensive kernel to compute. 
We expect a higher cost of the stokeslet as compared to the rotlet due
to the slower decay of the Fourier space part, as given
in table \ref{tab:est_fourier}. This means that larger FFT grids are 
needed to obtain the same accuracy. See e.g. the discussion in
connection to the left plot in figure \ref{fig:runtime_all_kernels}
where the choice of $M$ for the box $L=2$ is $48$ for the stokeslet
and $38$ for the rotlet. 

\begin{figure}[htbp]
  \centering
\includegraphics[width=0.49\textwidth,height=0.5\textwidth]{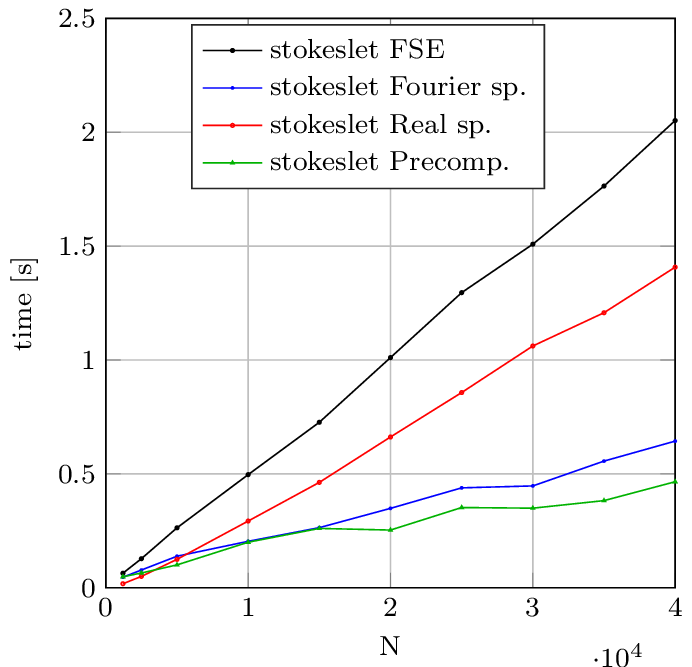}  
\includegraphics[width=0.49\textwidth,height=0.49\textwidth]{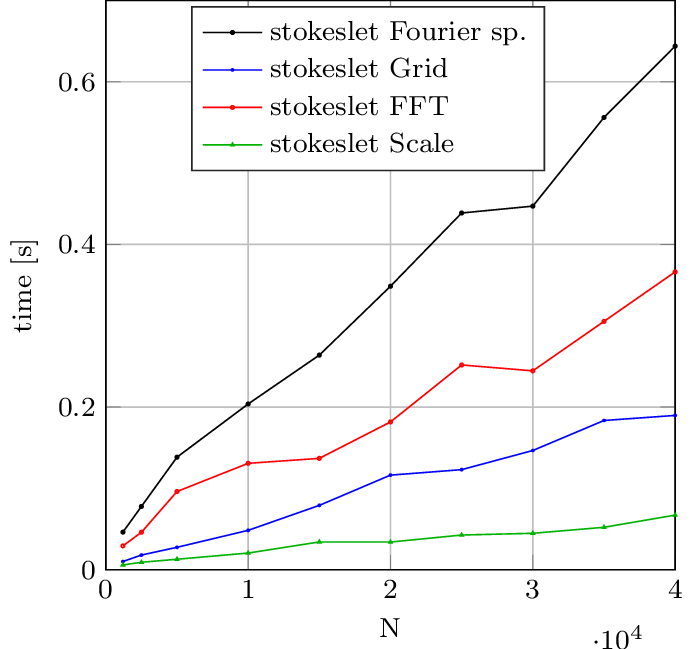}
  \caption{Breakdown of runtimes (left) and Fourier space runtime
    (right) for evaluating the Stokeslet as a function of number of
    particles. The full run-time was also shown in the left plot of
    figure \ref{fig:runtime_all_kernels}. }
\label{fig:stokeslet_detail}
\end{figure}

For the same system as in figure \ref{fig:runtime_all_kernels} (left), we now study the
computational cost for the different part of the calculations for the
stokeslet. In the left plot of figure \eqref{fig:stokeslet_detail} we
show the total evaluation runtime for the stokeslet sum together with
the three parts that makes up this total cost: the real space and
Fourier space evaluations plus the precomputation in Fourier space.
We use the choice of $\deltaL=d$ in \eqref{eq:deltaL_cond}, such that
$\Ltilde=L+d$. With this, $\Mtilde=M+P$, 
and the FFTs in the Fourier space evaluations will be of size
$(2\Mtilde)^3$. For the precomputation, the size of the FFT grids in
each dimension will be taken as the smallest even number that is
greater than $(1+\sqrt{3}) \Mtilde$.  The plot shows that the
computational cost is very similar for each part.  
As discussed above, the precomputation does not depend on the sources
and can be done only once as long as the domain size does not change.
Excluding the precomputation cost from the timing of the
stokeslet, the runtime is reduced somewhere between a quarter and one
third. Readjustment of $\xi$ to instead balance computational costs
excluding the precomputation, would yield a further reduced run time.

In the right plot of figure \ref{fig:stokeslet_detail}, we further
break down the cost of evaluating the  Fourier space sum into three
parts: {\em Grid} (the to and from grid operations with Gaussians),
{\em FFT} (the
total of $6$ FFTs) and {\em Scale}, the multiplication in step 3 of the
algorithm in section \ref{sec:discretization}. 
Note here that the oscillations in the FFT curve
are due to the fact that the FFT is more efficient for some grid sizes.
The scaling step is clearly the cheapest of the three parts.
The cost of the gridding step is $O(P^3N)$, where $P^3$ are the number
of grid points in the support of a Gaussian, and the cost of each
FFT of size $(2\Mtilde)^3$ is $O(\Mtilde^3 \log \Mtilde)$. 
Due to the connection to the real space sum, the choice of $\Mtilde$
will be such that this cost scales as $O(N \log N)$, as discussed in
section \ref{sec:comp-complexity}. 

\section{Conclusions}

We have presented a new fast summation method for free space Green's
functions of Stokes flow. The method is based on an Ewald
decomposition to split the sum in two parts, one in real space and
one in Fourier space. The real space sum can simply be truncated
outside of some radius of interaction that depends on the choice of
decomposition parameter and the required accuracy.  The focus of this
paper is on the Fourier space sum, the treatment of which is set in
the framework of the Spectral Ewald method, previously developed for
periodic problems \cite{Lindbo2010,AfKlinteberg2014a}.  The adaptation
to the free space problem involves a very recent approach to solving
the free-space harmonic and biharmonic equations using FFTs on a
uniform grid \cite{Vico2016}.  The Ewald Fourier space kernels for the
stokeslet, stresslet and rotlet are defined from the precomputed
Fourier representation of mollified harmonic (rotlet) and biharmonic
(stokeslet and stresslet) kernels, and the method can easily be
extended to any kernel that can be expressed as a differentiation of
the harmonic and/or biharmonic kernel.
New truncation error estimates have been derived for the free space
kernels. 

The extension of the FFT based Spectral Ewald method to the free space
problem incurs an additional computational cost compared to the
periodic problem. 
This is essentially due to the computation of larger FFTs, as
computational grids are zero-padded to the double size before the FFTs
are computed.  There is also an additional cost of two oversampled
FFTs for precomputing the Fourier representation of the mollified
harmonic or biharmonic kernel. This precomputation do not depend on
the sources, and the cost can often be amortized over many sum
evaluations.

Truncation error estimates have been derived for the kernels for
which they did not already exist, such that precise estimates of the
errors introduced by truncating the real and Fourier space sums are
available for all three kernels, the stokeslet, stresslet and rotlet. 
Errors decay exponentially in the physical distance and wave mode
number used for cut-off. 
Approximation errors in the evaluation of the Fourier sum decays
exponentially with the support of the Gaussians. 
An intricate detail needed to preserve the decoupling between
truncation and approximation errors that is not relevant for the periodic
Spectral Ewald method was discussed in section
\ref{sec:approximation-errors}. 

Numerical results are presented for the evaluation of the stokeslet,
stresslet and rotlet sums. They show the expected $O(N \log N)$ 
computational cost of the method. 
We have compared to an open source implementation of the FMM method
\cite{fmm_webpage}, and have shown that our method is competitive, as
it performs better for the uniform source distributions and high
accuracies considered here.  

With this, we have developed a new FFT based method for the fast
evaluation of free space Green’s functions for Stokes flow
(stokeslets, stresslets and rotlets) in a free space setting. This
free space Spectral Ewald method
allows the use of the same framework as the periodic one, which makes
it easy to swap methods depending on the problem under
consideration. The source code for the triply periodic SE method is
available online \cite{se_github}, and we plan to shortly release also the code for
this free space implementation. 

\section*{Acknowledgements}
  This work has been supported by the G\"{o}ran Gustafsson Foundation
  for Research in Natural Sciences and Medicine and by the Swedish
  Research Council under grant no. 2011-3178.  The authors gratefully
  acknowledge this support.


\bibliography{library}
\bibliographystyle{abbrvnat_mod}

\end{document}

%% file: defs.tex
\def\vec#1{\ensuremath{\mathchoice
                     {\mbox{\boldmath$\displaystyle\mathbf{#1}$}}
                     {\mbox{\boldmath$\textstyle\mathbf{#1}$}}
                     {\mbox{\boldmath$\scriptstyle\mathbf{#1}$}}
                     {\mbox{\boldmath$\scriptscriptstyle\mathbf{#1}$}}}}

\renewcommand{\v}[1]{\vec{#1}}
\newcommand{\expk}[1]{e^{-i\kb \cdot #1}}

\newtheorem{lemma}{Lemma}

\newcommand{\reals}{\mathbb R}
\newcommand{\Green}{G}
\newcommand{\stokeslet}{S}
\newcommand{\stresslet}{T}
\newcommand{\rotlet}{\Omega}
\newcommand{\harmonic}{H}
\newcommand{\biharmonic}{B}
\newcommand{\kmax}{k_{\infty}}
\newcommand{\trunc}{\mathcal{R}}

\newcommand{\xb}{\vec{x}}
\newcommand{\yb}{\vec{y}}
\newcommand{\ub}{\vec{u}}
\newcommand{\gb}{\vec{g}}
\newcommand{\fb}{\vec{f}}
\newcommand{\id}{\vec{I}}
\newcommand{\rb}{\vec{r}}
\newcommand{\pb}{\vec{p}}
\newcommand{\kb}{\vec{k}}
\newcommand{\mi}{\texttt{m}}
\renewcommand{\ni}{\texttt{n}}

\newcommand{\Z}{\mathbb{Z}}
\newcommand{\ftrans}{{\mathcal F}}
\newcommand{\rect}{\operatorname{rect}}
\newcommand{\ordo}{\mathcal{O}}
\newcommand{\erf}{\operatorname{erf}}
\newcommand{\erfc}{\operatorname{erfc}}
\newcommand{\K}{\operatorname{K}}
\newcommand{\splitfun}{\Phi}
\newcommand{\domain}{\mathcal{D}}
\newcommand{\domaintilde}{\tilde{\mathcal{D}}}
\newcommand{\osamp}{{s_f}}
\newcommand{\supp}{d} 
\newcommand{\deltaL}{\delta_L} 

\newcommand{\Ltilde}{\tilde{L}}
\newcommand{\Mtilde}{\tilde{M}}
\newcommand{\Mupgrid}{M_g}
\newcommand{\gridfcn}{\gb}
\newcommand{\gridfcnsc}{\vec{w}}
\newcommand{\Gwoexp}{A}
\newcommand{\kbar}{\bar{k}}
